\documentclass{amsart}

\newtheorem{THM}{Theorem}[section]
\newtheorem{LMA}[THM]{Lemma}
\newtheorem{PROP}[THM]{Proposition}
\newtheorem{CORO}[THM]{Corollary}

\newtheorem{PROB}[THM]{Problem}

\numberwithin{equation}{section}

\usepackage{amsmath}
\usepackage{amssymb}
\usepackage{euscript}

\newcommand{\showon}{\begin{eqnarray}}
\newcommand{\showoff}{\end{eqnarray}}
\newcommand{\showonq}{\begin{eqnarray*}}
\newcommand{\showoffq}{\end{eqnarray*}}

\newcommand{\NN}{\mathbb{N}}
\newcommand{\ZZ}{\mathbb{Z}}
\newcommand{\QQ}{\mathbb{Q}}
\newcommand{\RR}{\mathbb{R}}
\newcommand{\CC}{\mathbb{C}}
\newcommand{\EE}{\mathcal{E}}
\newcommand{\drop}{\smallsetminus}
\newcommand{\goesto}{\rightarrow}
\newcommand{\lraw}{\longrightarrow}
\newcommand{\none}{\varnothing}
\newcommand{\Po}{\EuScript{P}}
\newcommand{\orn}{\omega}
\newcommand{\Orn}{\Omega}
\newcommand{\eps}{\varepsilon}
\newcommand{\RC}{\EuScript{R}}
\newcommand{\rev}{\overline}
\newcommand{\und}{\underline}
\newcommand{\Hom}{\mathrm{Hom}}
\newcommand{\rank}{\mathrm{rank}}
\newcommand{\coker}{\mathrm{coker}}
\newcommand{\supp}{\mathrm{supp}}
\newcommand{\cha}{\mathrm{char}}
\newcommand{\la}{\langle}
\newcommand{\ra}{\rangle}
\newcommand{\Flow}{\Phi}
\newcommand{\thfn}{\boldsymbol{\vartheta}}
\newcommand{\thpsi}{\boldsymbol{\psi}}
\newcommand{\ALPHA}{\boldsymbol{\alpha}}
\newcommand{\BETA}{\boldsymbol{\beta}}
\newcommand{\jj}{\mathbf{j}}
\newcommand{\xx}{\mathbf{x}}
\newcommand{\yy}{\mathbf{y}}

\newcommand{\cutedge}{\bullet\!\!-\!\!\bullet}
\newcommand{\loopgph}{\bullet\!\!\mathsf{O}}
\newcommand{\Cyc}{\mathsf{C}}

\newcommand{\tv}{\mathtt{tv}}
\newcommand{\cplx}{\boldsymbol{\varkappa}}

\renewcommand{\span}{\mathrm{span}}
\newcommand{\righthookarrow}{\hookrightarrow}

\begin{document}

\title[Flows in graphs]{The algebra of flows in graphs}
\author{David G. Wagner}
\address{Department of Combinatorics and Optimization\\
University of Waterloo\\
Waterloo, Ontario, Canada\ \ N2L 3G1}
\email{\texttt{dgwagner@math.uwaterloo.ca}}
\thanks{Research supported by the Natural
Sciences and Engineering Research Council of Canada under
operating grant OGP0105392.}

\begin{abstract}
We define a contravariant functor $K^\cdot$ from the category of finite
graphs and graph morphisms to the category of finitely generated graded
abelian groups and homomorphisms.  For a graph $X$, an abelian group
$B$, and a nonnegative integer $j$, an element of $\Hom(K^j(X),B)$ is a 
coherent family of $B$-valued flows on the set of all graphs obtained by
contracting some $(j-1)$-set of edges of $X$;  in particular, $\Hom(K^1(X),
\RR)$ is the familiar (real) ``cycle-space'' of $X$.  We show that
$K^\cdot(X)$ is torsion-free and that its Poincar\'e polynomial is
the specialization $t^{n-k}T_X(1/t,1+t)$ of the Tutte polynomial
of $X$ (here $X$ has $n$ vertices and $k$ components).  Functoriality
of $K^\cdot$ induces a functorial coalgebra structure on $K^\cdot(X)$;
dualizing, for any ring $B$ we obtain a functorial
$B$-algebra structure on $\Hom(K^\cdot(X),B)$.  When $B$ is commutative we
present this algebra as a quotient of a divided power algebra, leading to
some interesting inequalities on the coefficients of the above Poincar\'e
polynomial.  We also provide a formula for the theta
function of the lattice of integer-valued flows in $X$, and conclude with
ten open problems. 
\end{abstract}
\maketitle

\section{Introduction.}

Two ideas for defining algebraic invariants of graphs have been
particularly successful up to now:  these are the spectral and the
$K$-theoretic approaches.  The spectral theory begins by associating a
(usually Hermitian) matrix $M$ with a graph $X$, and proceeds by relating 
combinatorial structure in $X$ to the spectral decomposition of $M$
\cite{Bi,BCN,CDS,Go}.
The $K$-theoretic approach is the theory of the Tutte polynomial
and its many interesting specializations \cite{Bi,BrOx,Tu,Tu2,W2}.
A third method which seems
promising but has been relatively neglected is the categorical approach.
This idea is to define a functor from the category of graphs and graph
morphisms into some (algebraic, geometric, topological,...) category, and
to use these other structures to analyze the category of graphs.
The most notable example of this approach so far is Walker's 
functorial setting \cite{Wa} for Lov\'asz's proof \cite{Lo} of the
Kneser conjecture; see also Bj\"orner's survey \cite{Bj}.

Here we define such a functor $K^\cdot$ into the category of finitely
generated graded abelian groups and homomorphisms.  The definition
of $K^\cdot(X)$ is a formalization of Kirchhoff's First Law on $X$ and all
of its images under repeated contractions of edges; for this reason
we refer to $K^\cdot(X)$ as the ``Kirchhoff group'' of $X$.
Imagine that the edges of $X$ represent pipes, all of the same cross-sectional
area, which are full of water (an incompressible fluid) and connected at
the vertices.  A flow of water in this system is represented by assigning
a (real) velocity parallel to each edge; Kirchhoff First Law states the linear
equations on these velocities which require that mass be conserved at
each vertex.  The Kirchhoff group is constructed so that the vector space
of all real-valued flows in $X$ turns out to be $\Hom(K^1(X),\RR)$.
For any abelian group $B$ and any nonnegative integer $j$, the elements
of $\Hom(K^j(X),B)$ can be interpreted as coherent families of $B$-valued
flows on the set of images of $X$ after contracting each $(j-1)$-set of edges,
as explained in Section 1.  The number $d_j(X):=\rank\:K^j(X)$ is thus the
dimension of the vector space of such coherent families of real-valued flows.

The Kirchhoff group turns out to be a remarkable invariant of graphs,
and its algebraic properties impose several inequalities on the
numbers $d_j(X)$.  This is quite interesting since, as we see in
Theorem $3.1$, if $X$ is a graph with $n$ vertices, $m$ edges, and $k$
connected components then
\showon
d_0+d_1t+\cdots+d_{m}t^{m}=
t^{n-k}T_X\left(\frac{1}{t},1+t\right),
\showoff
in which $T_X(x,y)$ is the Tutte polynomial of $X$.  To be specific about
these inequalities, for positive integers $a$ and $j$
there is a unique expression 
$$a=\binom{a_j}{j}+\binom{a_{j-1}}{j-1}+\cdots+\binom{a_i}{i}$$
such that $a_j>a_{j-1}>\cdots>a_i\geq i>0$.  The
\emph{$j$-th pseudopower of $a$} is 
$$\psi_j(a):=\binom{a_j+1}{j+1}+\binom{a_{j-1}+1}{j}+\cdots+
\binom{a_i+1}{i+1},$$
and we also put $\psi_j(0):=0$.  In Corollary $4.6$ we show that if
$X$ has $k$ connected components and $\ell$ cut-edges then
\showon
\begin{array}{l}
d_0=1,\ d_1=m-n+k,\ \mathrm{and}\ d_{m-\ell}=1,\\
d_j\neq 0\ \mathrm{if\ and\ only\ if}\ 0\leq j\leq m-\ell,\\
\mathrm{and}\ \mathrm{if}\  1\leq j\leq m-\ell-1 \ \mathrm{then}\
0<d_{j+1}\leq\psi_j(d_j).
\end{array}
\showoff
In Corollary $4.7$ we show that if $Y$ is a maximal forest in $X$ and
the fundamental cycles of $X$ with respect to $Y$ have lengths
$r_1,...,r_d$, then for each $0\leq j\leq m-\ell$, 
\showon
d_j\leq[t^j]\prod_{i=1}^d\left(1+t+t^2+\cdots+t^{r_i}\right),
\showoff
with the notation on the right side indicating the
coefficient of $t^j$ in the polynomial shown.
In Corollary $4.9$ we show that
\showon
d_j=\binom{d_1+j-1}{j} \ \mathrm{for}\ \mathrm{all}\ 0\leq j\leq g(X),
\showoff
in which $g(X)$ denotes the girth of $X$.  In Corollary $4.11$ we show that
\showon
\begin{array}{l}
d_0\leq d_1\leq\cdots\leq d_{\lfloor (m-\ell)/2\rfloor},\\
\mathrm{and\ if}\ 0\leq j\leq (m-\ell)/2\ \mathrm{then}\  d_j\leq d_{m-\ell-j}.
\end{array}
\showoff
All of these statements are just coarse numerical consequences of 
the algebraic structure of Kirchhoff groups, the description of which
is the main purpose of this paper.

In Section $1$ we formalize the definitions we require from graph theory
and define the Kirchhoff group $K^\cdot(X)$.  This is somewhat involved,
in that we must orient each edge of $X$ arbitrarily in order to write
down Kirchhoff's First Law and then, to obtain a group which is independent
of the orientation, symmetrize the construction over all orientations.
As a warm-up exercise we show in Theorem $1.2$ that if $X$ and $Y$ have
at most one vertex in common then $K^\cdot(X\cup Y)= K^\cdot(X)\otimes
K^\cdot(Y)$. Theorem $1.3$ is a technical but essential fact.

In Theorem $2.1$ and the following remarks we show that the Kirchhoff
group is a contravariant functor.  By Proposition $2.2$, every graph
morphism $f:Y\goesto X$ can be factored as a composition of morphisms
each of which has one of four very simple forms.  Propositions $2.3$ to $2.6$
describe the structure of the group homomorphism
$f^*:K^\cdot(X)\goesto K^\cdot(Y)$ in each of these cases.

Theorem $3.1$ uses the foregoing to show that Kirchhoff groups are
torsion-free, and relates them to the Tutte polynomial as in $(0.1)$.
Since $K^\cdot(X)$ is torsion-free, the ranks $d_j(X)$ determine
$K^\cdot(X)$ up to isomorphism.  However, functoriality
of $K^\cdot$ induces a functorial coalgebra structure on $K^\cdot(X)$,
and the situation becomes much more interesting.

For any ring $B$, the coalgebra structure on $K^\cdot(X)$ induces a
$B$-algebra structure on $\Hom(K^\cdot(X),B)$, and this is functorial.
We use the notation $\Flow_\cdot(X,B)$ for $\Hom(K^\cdot(X),B)$,
and refer to it as the
``$B$-circulation algebra of $X$''.  For commutative $B$, Theorem
$4.4$ gives a spanning set for $\Flow_\cdot(X,B)$ as a $B$-module;
if $B$ is also a $\QQ$-algebra then Corollary $4.5$ identifies a subset
of $\Flow_1(X,B)$ which generates $\Flow_\cdot(X,B)$ as a $B$-algebra.
This implies most of $(0.2)$ above (some parts require separate arguments).
The inequalities $(0.3)$ also follow, as does a sharpening of $(0.2)$ by
applying the Clements-Lindstr\"om Theorem.  For a commutative ring $B$,
Theorem $4.8$ presents $\Flow_\cdot(X,B)$ as a quotient of a divided
power algebra, and the equalities $(0.4)$ follow from the form of this
presentation.  Theorem $4.10$ describes a situation in which
there are monomorphisms $\Flow_j(X,B)\goesto\Flow_{m-\ell-j}(X,B)$ for
$0\leq j\leq (m-\ell)/2$ which, if $B$ is a commutative $\QQ$-algebra,
factor through all the intermediate graded components of $\Phi_\cdot(X,B)$.
This implies the inequalities $(0.5)$ above.

There is a natural Euclidean bilinear form $\la\cdot,\cdot\ra$
on $\Flow_\cdot(X,\RR)$, and it is integer-valued when restricted to
$\Flow_\cdot(X,\ZZ)$.  The theta function of the lattice
$\Flow_j(X,\ZZ)$ is the generating function for coherent families
$\varphi$ of integer-valued flows on the set of $(j-1)$-fold
contractions of $X$ with respect to the squared norm
$\la\varphi,\varphi\ra$.  Theorem $5.4$ gives a formula
for this theta function in the simplest nontrivial case, for
$\Flow_1(X,\ZZ)$.

In Section $6$ we conclude with ten open problems of varying feasibility
and interest.

\begin{center}
\textsc{Acknowledgments}
\end{center}

This research was performed in the spring of 1997, while I visited the
University of Notre Dame on sabbatical leave from the University of
Waterloo.  I am grateful to all the people at both institutions who helped
to make this possible.  After hearing me report on some related work
on partial orders at Oberwolfach in January 1996, Sergey Fomin suggested
that the relevant definitions made sense for graphs, a remark which
hinted in the direction of this paper and which I sincerely appreciate.

\section{The Kirchhoff group of a graph.}

The Kirchhoff group of a graph is a universal test object.  Informally,
it is a collection of ammeters, one for each edge of each contraction
of $X$.  The edges of the contractions of $X$ are identified with one another
in a natural way, so the readings on these ammeters must agree.
Kirchhoff First Law specifies further linear equations which must be satisfied
by the readings on these ammeters if they are measuring a collection of
flows in the contractions of $X$; these are the relations which are
built into the definition.

Let $X=(V,A,o,t,\rev{\cdot})$ be a general finite undirected graph,
by which we mean a finite set $V$ of \emph{vertices}, a finite set $A$
of \emph{arcs} (directed edges), an \emph{origin function} $o:A\goesto V$, a
\emph{terminus function} $t:A\goesto V$, and an \emph{arc reversal function}
$\rev{\cdot}:A\goesto A$, subject to the conditions that for all $a\in A$,
$\rev{a}\neq a$ and $\rev{\rev{a}}=a$ and $o(\rev{a})=t(a)$.
It follows that for all $a\in A$, $t(\rev{a})=o(a)$.  The set $E$ of
\emph{edges} of $X$ is the set of orbits of $\rev{\cdot}$ acting on $A$.
For $a\in A$ let $\und{a}:=\{a,\rev{a}\}$ be the corresponding edge.
For $v\in V$ let $X(v):=\{a\in A:\ t(a)=v\}$.
An \emph{orientation} of $X$ is a function $\orn:E\goesto A$ such that
$\orn(e)\in e$ for each $e\in E$;  the set of all orientations
of $X$ is denoted by $\Orn$.  For $\orn\in\Orn$ and $a\in A$, the
\emph{sign of $a$ relative to $\orn$} is
$$s_\orn(a):=\left\{\begin{array}{rl}
1  & \mathrm{if}\ \orn(\und{a})=a,\\
-1 & \mathrm{if}\ \orn(\und{a})=\rev{a}.\end{array}\right.$$

For $e\in E$, the \emph{contraction of $e$ in $X$}
is the graph $X_e=(V',A',o',t',\tilde{\cdot})$ defined as follows.
Choose $a\in e$ (the choice does not matter).
Let $V':=(V\drop\{o(a),t(a)\})\cup\{e\}$, and define $\pi:V\goesto V'$ by
$\pi(v):=v$ if $v\not\in\{o(a),t(a)\}$, and $\pi(v):=e$ otherwise.  Let
$A':=A\drop\{a,\rev{a}\}$, let $o':=\pi\circ o|_{A'}$ and $t':=\pi\circ
t|_{A'}$, and let $\tilde{\cdot}:=\rev{\cdot}|_{A'}$.   Any $\orn\in\Orn$
induces an orientation on $X_e$ (by restriction); somewhat improperly
we shall also denote this orientation of $X_e$ by $\omega$.
For distinct edges $c,e\in E$ one sees that $(X_e)_c$ is
naturally isomorphic with $(X_c)_e$, and it follows that for
every subset $\sigma\subseteq E$, the contracted graph $X_\sigma$ is
well-defined.  (It may be helpful to think of a graph $X$ as a finite
CW-complex of dimension at most one, in which case $X_\sigma$ is the
quotient complex of $X$ modulo the subcomplex $Z$ with vertices $V$ and 
edges $\sigma$;  in particular, $V(X_\sigma)$ is the set of connected
components of this subgraph $Z$ of $X$.)

The set of all $j$-element subsets of $E$ is denoted by $\Po_j(X)$,
and $\Po(X):=\bigcup_j\Po_j(X)$.  Ordered by inclusion, $\Po(X)$ is
of course a Boolean lattice. For each integer $j$, let $F^j(X)$ be the
free abelian group with $\{X_\sigma:\ \sigma\in\Po_j(X)\}$ as a basis:
$$F^j(X) := \ZZ\{X_\sigma:\ \sigma\in\Po_j(X)\}.$$
Then $F^\cdot(X):=\bigoplus_j F^j(X)$ is a graded finitely generated
free abelian group.  Homomorphisms between graded groups are required
to preserve the grading.  For a graded abelian group $M^\cdot$ and
integer $r$,  $M^\cdot(r)$ denotes $M^\cdot$ with the grading shifted by
$r$:  the $j$-th graded component of $M^\cdot(r)$ is the $(j+r)$-th graded
component of $M^\cdot$.
Frequently, we define a homomorphism on a free module over a ring $B$
by defining its values on a basis only;  it is to be understood that the
definition is implicitly extended $B$-linearly to the whole module.  (The
case $B=\ZZ$ is that of free abelian groups.)

For each $\sigma\in\Po_j(X)$ there is a natural
inclusion
\showon
\hat{\eta}_{X,\sigma}:\ F^\cdot(X_\sigma)(-j) \righthookarrow F^\cdot(X)
\showoff
induced by the identifications $(X_\sigma)_\tau = X_{\sigma\cup\tau}$
for all $\tau\in\Po(X_\sigma)$.  If $\sigma,\rho\in\Po(X)$ are disjoint
then clearly
\showon
\hat{\eta}_{X,\sigma}\circ\hat{\eta}_{X_\sigma,\rho}
=\hat{\eta}_{X,\sigma\cup\rho}.
\showoff

Fix an orientation $\orn\in\Orn$.  For each vertex $v\in V$ define a
relation
$$R_\orn(X,v) := \sum_{a\in X(v)} s_\orn(a) X_{\und{a}}$$
in $F^\cdot(X)$, and let $\RC_\orn(X):=\{R_\orn(X,v):\ v\in V(X)\}$. 
These are Kirchhoff's Laws on $X$.
For each integer $j$ define a subgroup of $F^j(X)$ by
$$N_\orn^j(X) := \span_\ZZ\bigcup_{\sigma\in\Po_{j-1}(X)}
\hat{\eta}_{X,\sigma}(\RC_\orn(X_\sigma)).$$
Then $N_\orn^\cdot(X):=\bigoplus_j N_\orn^j(X)$ is a graded subgroup of
$F^\cdot(X)$.  The \emph{Kirchhoff group of $X$ relative to $\orn$} is
the quotient
\showon
K_\orn^\cdot(X) &:=& F^\cdot(X)/N_\orn^\cdot(X).
\showoff

To define a group $K^\cdot(X)$ which depends naturally on $X$,
and not on any particular orientation of $X$, we symmetrize the above
construction using automorphisms $\psi_\eps^\orn:\ F^\cdot(X)\goesto
F^\cdot(X)$ for all $\orn,\eps\in \Orn$, defined as follows.
For each $\sigma\in\Po(X)$, let
$m_\sigma(\orn,\eps):=\#\{e\in\sigma: \orn(e)\neq\eps(e)\}$, and
define $\psi_\eps^\orn(X_\sigma):= (-1)^{m_\sigma(\orn,\eps)} X_\sigma$.
Let $\Psi:=\{\psi_\eps^\orn:\ \orn,\eps\in\Orn\}$.  (In fact, $\Psi$ is
an elementary abelian group of order $2^{\#E}$, but we do not use this
fact.)  One easily sees that for all $\orn\in\Orn$, $\psi_\orn^\orn=1$,
and that for all $\orn,\eps,\nu\in\Orn$, $\psi_\nu^\eps\circ
\psi_\eps^\orn= \psi_\nu^\orn$.  Let $F^\cdot_\Orn(X)$ be the subgroup
of $\bigoplus_{\orn\in\Orn} F^\cdot(X)$ consisting of all elements
$(T_\orn:\ \orn\in\Orn)$ such that $T_\eps=\psi_\eps^\orn(T_\orn)$
for all $\orn,\eps \in\Orn$.  For each $\orn\in\Orn$ there is an
isomorphism
\showon
\hat{\kappa}_\orn:\ F^\cdot(X)\goesto F_\Orn^\cdot(X)
\showoff
defined by $\hat{\kappa}_\orn(T):= (\psi_\eps^\orn(T):\ \eps\in\Orn)$
for all $T\in F^\cdot(X)$.  These are compatible with $\Psi$
in the sense that for any $\orn,\eps\in\Orn$, $\hat{\kappa}_\orn=
\hat{\kappa}_\eps\circ\psi^\orn_\eps$.
For any $\orn,\eps\in\Orn$, any $\sigma\in\Po(X)$, and
any $v\in V(X_\sigma)$, one sees that $\psi_\eps^\orn(R_\orn(X_\sigma,v))
=(-1)^{m_\sigma(\orn,\eps)}R_\eps(X_\sigma,v)$.  From this it follows
that for any $\orn,\eps\in\Orn$, $\psi_\eps^\orn(N_\orn^\cdot(X))=
N_\eps^\cdot(X)$.  Therefore, the subgroup $N_\Orn^\cdot(X):=
\hat{\kappa}_\orn(N_\orn^\cdot(X))$ of $F_\Orn^\cdot(X)$ is independent
of the choice of $\orn\in\Orn$.  Finally, we define \emph{the Kirchhoff
group of $X$} to be the quotient
$$K^\cdot(X):=F_\Orn^\cdot(X)/N_\Orn^\cdot(X).$$
From $(1.4)$ we obtain a family of isomorphisms
\showon
\kappa_\orn:\ K_\orn^\cdot(X)\goesto K^\cdot(X)
\showoff
for $\orn\in\Orn$ which are compatible with $\Psi$.
If $\mu_\orn:K_\orn^\cdot(X)\goesto M^\cdot$ is a family of homomorphisms
indexed by $\Orn$ and compatible with $\Psi$, then the homomorphism
$\mu:=\mu_\orn\circ\kappa_\orn^{-1}$ is independent of the choice of
$\orn\in\Orn$.  By virtue of this, in any particular argument one is
free to choose $\orn\in\Orn$ arbitrarily, and to work with
$K_\orn^\cdot(X)$ instead of $K^\cdot(X)$, as long as one checks
that the construction is compatible with $\Psi$.

A word of explanation is in order as to what the Kirchhoff group really
means.  For each $\sigma\in\Po(X)$ and $e\not\in\sigma$, think of
$X_{\sigma\cup\{e\}}$ as representing the edge $e$ in $X_\sigma$.  Each
$X_\tau$ thus represents one edge in each of $\#\tau$ different contractions
of $X$.  Fixing a reference orientation $\orn\in\Orn$, the set
$\RC_\orn(X_\sigma)$ is a formulation of Kirchhoff First Law on the
graph $X_\sigma$.
Thus, for any integer $j$ and abelian group $B$, an element of
$\Hom(K^j(X),B)$ is a coherent family of $B$-valued flows on
$\{X_\sigma:\ \sigma\in\Po_{j-1}(X)\}$.  This $\Hom(K^j(X),B)$ is naturally
an abelian group, and moreover if $B$ is a ring then it is naturally a
$B$-module.  In particular, $\Hom(K^1(X),\RR)$ is the vector space of
real-valued flows on $X$; that is, the (real) \emph{cycle space} of $X$.

The following very simple lemma is a key ingredient of many proofs in
this paper.  For a proposition $P$ we define the \emph{truth
value of $P$} to be $\tv[P]:=1$ if $P$ is true and $\tv[P]:=0$ if $P$ is
false.
\begin{LMA}  Let $X$ be a graph.  For any $\orn\in\Orn$ and
$U\subseteq V$,
$$\sum_{v\in U} R_\orn(X,v) = \sum_{a\in X(U)} s_\orn(a) X_{\und{a}},$$
in which $X(U):=\{a\in A:\ t(a)\in U\ \mathrm{and}\ o(a)\not\in U\}$.
\end{LMA}
\begin{proof}
For each $e\in E$, $X_e$ appears on the left side with coefficient
$\tv[t(\orn(e))\in U]-\tv[o(\orn(e))\in U]$, from which the
result follows.
\end{proof}

For any graph $X$, the assignment $X_\none\mapsto 1$ defines a natural
isomorphism $F^0(X)=\ZZ$.  Since $X_\none$ is fixed by every element
of $\Psi$, its image $u_X\in K^0(X)$ is well-defined.
This gives a natural isomorphism $K^0(X)=\ZZ$ by sending $u_X$ to $1$.
If $X$ and $Y$ are graphs with disjoint edge-sets, then the bilinear
map $F^\cdot(X)\times F^\cdot(Y)\goesto F^\cdot(X\cup Y)$ defined by
$(X_\sigma,Y_\tau)\mapsto(X\cup Y)_{\sigma\cup\tau}$ for all
$\sigma\in\Po(X)$ and $\tau\in\Po(Y)$ induces a natural
isomorphism $F^\cdot(X)\otimes F^\cdot(Y)=F^\cdot(X\cup Y)$.
(Tensor products are over $\ZZ$ unless otherwise noted.)

\begin{THM}  Let $X$ and $Y$ be graphs with disjoint edge-sets, and with
at most one vertex in common.  Then $K^\cdot(X\cup Y)=
K^\cdot(X)\otimes K^\cdot(Y).$ \end{THM}
\begin{proof}
Fix reference orientations $\orn\in\Orn(X)$ and $\gamma\in\Orn(Y)$;
these determine an orientation $\nu\in\Orn(X\cup Y)$.
Consider a relation $R_\orn(X_\sigma,v)$ (so $\sigma\in\Po(X)$ and
$v\in V(X_\sigma)$) for which $v$ is not the image in $V(X_\sigma)$
of any vertex in $V(Y)$, and let $\tau\in\Po(Y)$.  Then
\showonq
R_\orn(X_\sigma,v)\otimes Y_\tau &=&
\sum_{a\in X_\sigma(v)} s_\orn(a) (X\cup Y)_{\sigma\cup\{\und{a}\}\cup\tau}\\
&=& \sum_{a\in(X\cup Y)_{\sigma\cup\tau}(v)} s_\orn(a)
(X\cup Y)_{\sigma\cup\tau\cup\{\und{a}\}}\\
&=&R_{\nu}((X\cup Y)_{\sigma\cup\tau},v),
\showoffq
which is in $N_\nu^\cdot(X\cup Y)$.  But there is at most one vertex
$w\in V(X_\sigma)$ which is the image of some vertex of $V(Y)$.  By
Lemma $1.1$, $\sum_{v\in V(X_\sigma)} R_\orn(X_\sigma,v)=0$, and so if
such a $w$ exists then $R_\orn(X_\sigma,w)\otimes Y_\tau$ is also in
$N_\nu^\cdot(X\cup Y)$.  By the symmetry exchanging $X$ and $Y$ we
conclude that $N_\orn^\cdot(X)\otimes F^\cdot(Y)+F^\cdot(X)\otimes
N_\gamma^\cdot(Y)\subseteq N_\nu^\cdot(X\cup Y)$, and so we have
an epimorphism 
$$\phi:K_\orn^\cdot(X)\otimes K_\gamma^\cdot(Y)\goesto
K_\nu^\cdot(X\cup Y).$$
Now consider any $\rho\in\Po(X\cup Y)$, and let $\sigma:=\rho\cap E(X)$
and $\tau:=\rho\cap E(Y)$.  For $v$ a vertex of $(X \cup Y)_\rho$ which
is not in $V(X_\sigma)\cap V(Y_\tau)$ we may assume that $v\in V(X_\sigma)$,
by symmetry.  Then $R_\orn(X_\sigma,v) \otimes Y_\tau=
R_\nu((X \cup Y)_\rho,v)$.  But there is at most one vertex $w$ in
$V(X_\sigma)\cap V(Y_\tau)$, and by Lemma 1.1,
$\sum_{w\in V((X\cup Y)_\rho)} R_\nu((X\cup Y)_\rho,v)=0$.
This shows that $N_\nu^\cdot(X\cup Y)\subseteq
N_\orn^\cdot(X)\otimes F^\cdot(Y)+F^\cdot(X)\otimes N_\gamma^\cdot(Y)$,
and so $\phi$ is an isomorphism.  To complete the proof one must check that
the construction of $\phi$ is compatible with $\Psi(X)$ and $\Psi(Y)$ and
the isomorphisms $\kappa_\nu:K_\nu^\cdot(X\cup Y)\goesto K^\cdot(X\cup Y)$,
but this is straightforward.  \end{proof}

Now fix an $\orn\in\Orn$, and consider one of the relations
$R_\orn((X_\sigma)_\tau,v)$ of $N_\orn^\cdot(X_\sigma)$,
so $\tau\in\Po(X_\sigma)$ and $v\in V((X_\sigma)_\tau)$.  Notice that
$$\hat{\eta}_{X,\sigma}(R_\orn((X_\sigma)_\tau,v))=
\sum_{a\in (X_\sigma)_\tau(v)} s_\orn(a) X_{\sigma\cup\tau\cup\{\und{a}\}}
=R_\orn(X_{\sigma\cup\tau},v).$$
From this it follows that for all $\orn\in\Orn$ and $\sigma\in\Po_j(X)$,
$\hat{\eta}_{X,\sigma}(N_\orn^\cdot(X_\sigma)(-j))\subseteq
N_\orn^\cdot(X),$ and so $(1.1)$ induces a family of homomorphisms
\showon
\eta^\orn_{X,\sigma}:\ K_\orn^\cdot(X_\sigma)(-j)\goesto K_\orn^\cdot(X)
\showoff
which are compatible with $\Psi$ and the maps $\kappa_\orn$.
Hence, these induce a natural homomorphism 
$$ \eta_{X,\sigma}:K^\cdot(X_\sigma)(-j) \goesto K^\cdot(X).$$
Theorem $1.3$ has important consequences for the structure of Kirchhoff
groups and circulation algebras, as we see in Sections $3$ and $4$.
\begin{THM}  Let $X$ be a graph and let $\sigma\in\Po_j(X)$.\\
\textup{(a)}\ If $\sigma$ contains a cut-edge of $X$ then
$\eta_{X,\sigma}=0$.\\
\textup{(b)}\ If $\sigma$ does not contain any cut-edges of $X$ then
there is a homomorphism $\pi:K^\cdot(X)\goesto
K^\cdot(X_\sigma)(-j)$ such that $\pi\circ\eta_{X,\sigma}=1$.
In particular, $\eta_{X,\sigma}$ is injective.\end{THM}
\begin{proof}
Fix an orientation $\orn\in\Orn$.  Identify $F^\cdot(X_\sigma)(-j)$
with its image under $\hat{\eta}_{X,\sigma}$.
From $(1.1)$, $(1.6)$, and the kernel-cokernel exact sequence
(see, \emph{e.g.}, Lemma $II.5.2$ of Mac Lane \cite{Mac}) , we see that
$$0\lraw \ker(\eta^\orn_{X,\sigma})\lraw
\frac{N^\cdot_\orn(X)}{N^\cdot_\orn(X_\sigma)(-j)} \stackrel{\phi}{\lraw}
\frac{F^\cdot(X)}{F^\cdot(X_\sigma)(-j)} \lraw
\coker(\eta^\orn_{X,\sigma}) \lraw 0$$
is exact, and $\phi$ is induced by the inclusion $N_\orn^\cdot(X)
\righthookarrow F^\cdot(X)$.  Hence,
\showon
\ker(\eta^\orn_{X,\sigma})\simeq\ker(\phi)=
\frac{N^\cdot_\orn(X)\cap F^\cdot(X_\sigma)(-j)}
{N^\cdot_\orn(X_\sigma)(-j)}.
\showoff

For part (a), assume that $e\in\sigma$ is a cut-edge of $X$, and let
$\rho:=\sigma\drop\{e\}$.  Then $e$ is a cut-edge of $X_\tau$ for
every $e\not\in\tau\in\Po(X)$.  To prove that $\eta^\orn_{X,\sigma}=0$,
it suffices to show that $F^\cdot(X_\sigma)(-j)\subseteq N^\cdot_\orn(X)$.
Consider any $\tau\in\Po(X_\sigma)$, and the corresponding element
$X_{\sigma\cup\tau}$ of
$F^\cdot(X_\sigma)(-j)$.  Now $e$ is a cut-edge of $X_{\rho\cup\tau}$.
Let $U$ be the set of vertices of $X_{\rho\cup\tau}$ in the same
connected component of $X_{\rho\cup\tau}\drop\{e\}$  as $t(\orn(e))$.
By Lemma 1.1
we have $\sum_{v\in U} R_\orn(X_{\rho\cup\tau},v)=X_{\sigma\cup\tau}
\in N_\orn^\cdot(X)$.  Since $F^\cdot(X_\sigma)(-j)$ is generated by
elements of this form, $\eta^\orn_{X,\sigma}=0$.  The isomorphism
$\kappa_\orn$ shows that $\eta_{X,\sigma}=0$.

It suffices to prove part (b) in the special case $j=1$, because of
$(1.2)$.  Accordingly, $\sigma=\{e\}$ and $e$ is not a cut-edge of $X$.
Let $C\subseteq X$ be a cycle in $X$ which contains $e$, and let $\zeta:=
E(C)$.  We may choose $\orn\in\Orn$ so that for each $v\in V(C)$ there
is exactly one $c\in\zeta$ such that $v=t(\orn(c))$.
Define $\hat{\pi}:F^\cdot(X)\goesto F^\cdot(X_e)(-1)$ by 
$$\hat{\pi}(X_\tau):=\left\{\begin{array}{ll}
X_\tau & \mathrm{if}\ e\in\tau,\\
\sum_{c\in\tau\cap\zeta}X_{\tau\cup\{e\}\drop\{c\}} & \mathrm{if}\
e\not\in\tau,\end{array}\right.$$
for all $\tau\in\Po(X)$.  It is clear that $\hat{\pi}:F^\cdot(X)\goesto
F^\cdot(X_e)(-1)$ and that $\hat{\pi}\circ\hat{\eta}_{X,e}=1$.
To finish the proof it remains to show that
$\hat{\pi}(N_\orn^\cdot(X))\subseteq N_\orn^\cdot(X_e)(-1)$, since
then there is an induced homomorphism $\pi_\orn: K^\cdot_\orn(X)\goesto
K^\cdot_\orn(X_e)(-1)$.  Having done this, it is straightforward to
check compatibility with $\Psi$ to obtain $\pi$ as claimed.

Consider a relation $R_\orn(X_\tau,v)$ in $N_\orn^\cdot(X)$, so $\tau\in
\Po(X)$ and $v\in V(X_\tau)$.  If $e\in\tau$ then
$\hat{\pi}(R_\orn(X_\tau,v)) =R_\orn(X_\tau,v)$, as is clear.
If $e\not\in\tau$ then
\showonq
\hat{\pi}(R_\orn(X_\tau,v)) &=&
\left(\tv[t(e)=v]-\tv[o(e)=v]\right)X_{\tau\cup\{e\}}\\
& & +\sum_{a\in X_\tau(v)\drop e}
s_\orn(a)\sum_{c\in(\tau\cup\{\und{a}\})\cap\zeta}
X_{\tau\cup\{\und{a}\}\cup\{e\}\drop\{c\}}\\
&=&\sum_{c\in\tau\cap\zeta}\sum_{a\in X_\tau(v)} s_\orn(a)
X_{\tau\cup\{e\}\cup\{\und{a}\}\drop\{c\}}+\sum_{a\in X_\tau(v):\
\und{a}\in\zeta} s_\orn(a) X_{\tau\cup\{e\}}.
\showoffq
The second sum is zero, since it is either empty (if $v$ is not a
vertex of the image of $C$ in $X_\tau$) or it has an even number of
terms, and these cancel in pairs because the orientation $\orn$ is
directed consistently around $C$.  For $c\in\tau\cap\zeta$,
let $U(c)$ denote the set of vertices of $X_{\tau\cup\{e\}\drop\{c\}}$
which are contracted to the image of $v$ in $X_{\tau\cup\{e\}}$;
this $U(c)$ has either one or two elements.  By Lemma 1.1,
$$\sum_{a\in X_\tau(v)} s_\orn(a)X_{\tau\cup\{e\}\cup\{\und{a}\}\drop\{c\}}
=\sum_{u\in U(c)} R_\orn(X_{\tau\cup\{e\}\drop\{c\}},u),$$
and we conclude that $\hat{\pi}(R_\orn(X_\tau,v))\in N^\cdot_\orn(X)$,
as required.
\end{proof}
Notice that the homomorphism $\pi$ constructed in part (b) is not
natural, as it depends on arbitrary choices of cycles in $X$.

\section{Functorial properties of Kirchhoff groups.}

A \emph{graph morphism} $f:Y\goesto X$ is a pair of functions
$f_V:V(Y)\goesto V(X)$ and $f_A:A(Y)\goesto A(X)$ which are compatible
with the functions $o$, $t$, $\rev{\cdot}$ defining $Y$ and $X$.
These induce a corresponding function $f_E:E(Y)\goesto E(X)$ as well.
We construct, for each graph morphism $f:Y\goesto X$, an induced
homomorphism $f^*:K^\cdot(X)\goesto K^\cdot(Y)$, in such a way that
$K^\cdot$ becomes a contravariant functor from the category of finite
graphs and graph morphisms to the category of finitely generated
graded abelian groups and homomorphisms.   Then we investigate the
relationship between structural properties of $f$ and algebraic
properties of $f^*$.

For $\sigma\in\Po(X)$, the set $T_f(\sigma)$ of \emph{$f$-transversals
of $\sigma$} is the set of all $\tau\in\Po(Y)$ such that for each
$e\in\sigma$ there is exactly one $c\in\tau$ for which $f_E(c)=e$,
and for each $e\not\in\sigma$ there is no $c\in\tau$ for which $f_E(c)=e$.
In other words, $T_f(\sigma)$ may be naturally identified with
$\prod\{f^{-1}_E(e):\ e\in\sigma\}$.  Define $f^\sharp:F^\cdot(X)
\goesto F^\cdot(Y)$ by
\showon
f^\sharp(X_\sigma):=\sum_{\tau\in T_f(\sigma)}Y_\tau
\showoff
for each $\sigma\in\Po(X)$.
Given any $\orn\in\Orn(X)$, the \emph{pullback of $\orn$ along $f$}
is the unique $\gamma\in\Orn(Y)$ making the following diagram
commute:
\showonq
\begin{array}{ccc}
E(Y) & \stackrel{f_E}{\longrightarrow} & E(X)\\
\gamma\downarrow & & \downarrow\orn\\
A(Y) & \stackrel{f_A}{\longrightarrow} & A(X)\end{array}
\showoffq

\begin{THM}  Let $f:Y\goesto X$ be a graph morphism.  Fix $\orn\in\Orn(X)$
and let $\gamma$ be the pullback of $\orn$ along $f$.  Then
$f^\sharp(N_\orn^\cdot(X))$ is a subgroup of $N_\gamma^\cdot(Y)$.
\end{THM}
\begin{proof}
Consider one of the relations $R_\orn(X_\sigma,v)$ of $N_\orn^\cdot(X)$,
so $\sigma\in\Po(X)$ and $v\in V(X_\sigma)$.  Let $U'\subseteq V(X)$
be the set of vertices of $X$ which are contracted to $v$ in $X_\sigma$,
and let $U:=f_V^{-1}(U')\subseteq V(Y)$ be the preimage of $U'$ in $Y$.
For each $\rho\in T_f(\sigma)$ let $U(\rho)\subseteq V(Y_\rho)$ be the
set of vertices of $Y_\rho$ to which some vertex in $U$ is contracted.
Let $\zeta:=f_E^{-1}(\sigma)$ be the set of edges of $Y$ mapped to some
edge of $\sigma$.  Notice that if $\rho\in T_f(\sigma)$ then $\rho\subseteq
\zeta$.  Now
\showonq
f^\sharp(R_\orn(X_\sigma,v)) &=&
\sum_{a\in X_\sigma(v)} s_\orn(a)
\sum_{\tau\in T_f(\sigma\cup\{\und{a}\})} Y_\tau\\
&=& 
\sum_{\rho\in T_f(\sigma)} \sum_{a\in X_\sigma(v)} s_\orn(a)
\sum_{c\in f_E^{-1}(\und{a})} Y_{\rho\cup\{c\}}.
\showoffq
As a subset of $A(X)$, $X_\sigma(v)$ is the set of arcs $a$ such that
$\und{a}\not\in\sigma$ and $t(a)\in U'$.   Thus, $f_A^{-1}(X_\sigma(v))$ 
is the set of arcs $a\in A(Y)$ such that $\und{a}\not\in\zeta$
and $t(a)\in U$.  Fix any $\rho\in T_f(\sigma)$, and let $A(\rho)$
denote the image of $f_A^{-1}(X_\sigma(v))$ in $Y_\rho$.  We have
$$\sum_{a\in X_\sigma(v)} s_\orn(a)\sum_{c\in f_E^{-1}(\und{a})}
Y_{\rho\cup\{c\}}=\sum_{a\in A(\rho)}s_\gamma(a) Y_{\rho\cup\{\und{a}\}}.$$
But $A(\rho)$ is the set of arcs $a\in A(Y_\rho)$ such that
$\und{a}\not\in\zeta$ and $t(a)\in U(\rho)$. If $\und{a}\in\zeta$ then
both $t(a)\in U(\rho)$ and $t(\rev{a})\in U(\rho)$.  Since
$s_\gamma(a)=-s_\gamma(\rev{a})$, the inclusion of such arcs into
the summation will cancel to zero.  Thus,  with $Y_\rho(U(\rho))$
denoting the set of arcs $a\in A(Y_\rho)$ such that $t(a)\in U(\rho)$
and $o(a)\not\in U(\rho)$,
\showonq
\sum_{a\in A(\rho)} s_\gamma(a) Y_{\rho\cup\{\und{a}\}} &=&
\sum_{a\in Y_\rho(U(\rho))} s_\gamma(a) Y_{\rho\cup\{\und{a}\}}\\
&=& \sum_{u\in U(\rho)} R_\gamma(Y_\rho,u),
\showoffq
by Lemma 1.1.  Putting the pieces together, we have
$$f^\sharp(R_\orn(X_\sigma,v))=
\sum_{\rho\in T_f(\sigma)} \sum_{u\in U(\rho)} R_\gamma(Y_\rho,u),$$
showing that $f^\sharp$ takes an arbitrary generator
of $N^\cdot_\orn(X)$ into $N^\cdot_\gamma(Y)$.
\end{proof}

For a graph morphism $f:Y\goesto X$ and orientation $\orn\in\Orn(X)$,
let $\gamma\in\Orn(Y)$ be the pullback of $\orn$ along $f$.  By Theorem
$2.1$ there is a homomorphism
$$f^\orn:K_\orn^\cdot(X)\goesto K_\gamma^\cdot(Y)$$
which is well-defined by putting $f^\orn(T+N_\orn^\cdot(X)):=
f^\sharp(T)+N_\gamma^\cdot(Y)$ for all $T\in F^\cdot(X)$.  These
are compatible with $\Psi(X)$ and the maps $\kappa_\gamma$ for
$\gamma\in\Orn(Y)$, and thus they define a homomorphism
$$f^*:K^\cdot(X)\goesto K^\cdot(Y).$$
If $g:Z\goesto Y$ is another graph morphism then it is clear from the
construction that $(f\circ g)^\sharp = g^\sharp \circ f^\sharp$, and
it follows that $(f\circ g)^* = g^* \circ f^*$.  This establishes that
$K^\cdot$ is a contravariant functor, as claimed.

\begin{PROP}  Every graph morphism $f:Y\goesto X$ can be factored as
$$Y\stackrel{g}{\lraw} Z \stackrel{h}{\lraw} W \stackrel{i}{\lraw} X$$
such that:\ $g_V$ is surjective and $g_A$ is bijective;
$h_V$ is bijective and $h_A$ is surjective; $i_V$ and $i_A$ are both
injective.\end{PROP}
\begin{proof}
Let $W$ be the image $f(Y)$ of $Y$ in $X$ and let $i:W\goesto X$ be
the inclusion of $W$ as a subgraph of $X$.  There is a graph morphism
$p:Y\goesto W$ such that $f=i\circ p$, and both $p_V$ and $p_A$ are
surjective.  Let $Z$ be the graph defined by $V(Z):=V(W)$, $A(Z):=A(Y)$,
$o_Z:=p_V\circ o_Y$, $t_Z:= p_V\circ t_Y$, and with the same arc-reversal
as $Y$.  The morphism $g$ is defined by $g_V:=p_V$ and
$g_A:=\mathrm{id}_{A(Y)}$.  The morphism $h$ is defined by
$h_V:=\mathrm{id}_{V(W)}$ and $h_A:=p_A$.
\end{proof}
(In fact, this factorization is unique up to isomorphism, but we make
no use of this fact.)

An \emph{elementary injection} is an injective function
$f:U\goesto W$ between finite sets such that $\#W=1+\#U$.
An \emph{elementary surjection} is a surjective function
$f:U\goesto W$ between finite sets such that $\#U=1+\#W$.
Let $f:Y\goesto X$ be a graph morphism, factored $f=i\circ h\circ g$
as in Proposition $2.2$.  We may further factor $g$ into a sequence of
morphisms $g'$ for each of which $g'_V$ is an elementary surjection
and $g'_A$ is bijective.  We may also factor $h$ into a sequence of
morphisms $h'$ for each of which $h'_V$ is bijective and $h'_E$ is
an elementary surjection.  Finally, we may factor $i$ into a sequence
of morphisms $i'$ for each of which exactly one of $i'_V$ and $i'_E$
is an elementary injection, and the other is a bijection.
By functoriality of $K^\cdot$, it suffices (in principle) to
describe $f^*$ for graph morphisms $f$ having one of these
special forms.

\begin{PROP} Let $i:Y\goesto X$ be a graph morphism such that
$i_V$ is an elementary injection and $i_A$ is bijective.  Then
$i^*$ is an isomorphism. \end{PROP}
\begin{proof}
The graph $X$ is the disjoint union of $Y$ and a graph $\bullet$ with
one vertex and no edges.  By Theorem $1.2$, $K^\cdot(X)=K^\cdot(Y)
\otimes K^\cdot(\bullet)= K^\cdot(Y)\otimes\ZZ=K^\cdot(Y)$.
This isomorphism is induced by $i_E$, so it is $i^*$.  \end{proof}

Proposition $2.4$ is an algebraic analogue of the deletion/contraction
algorithm of graph theory, and is a key component of several inductive
proofs in what follows.  In the following propositions it is convenient
to use the notation $f^\flat$ for the restriction of $f^\sharp$ to
$N^\cdot_\orn(X)$.

\begin{PROP} Let $i:Y\goesto X$ be a graph morphism such that
$i_V$ is bijective and $i_E$ is an elementary injection.  Let
$e$ be the unique edge of $X$ not in $i_E(E(Y))$.  Then
$$K^\cdot(X_e)(-1)\stackrel{\eta_{X,e}}{\lraw}
K^\cdot(X)\stackrel{i^*}{\lraw} K^\cdot(Y)\lraw 0$$
is exact.  \end{PROP}
\begin{proof}
Identify $V(X)$ and $V(Y)$ via $i_V$, and consider $E(Y)$ as a subset of
$E(X)$ via $i_E$.  For all $\sigma\in\Po(X)$,
$i^\sharp(X_\sigma):=Y_\sigma$ if $e\not\in\sigma$, and
$i^\sharp(X_\sigma):=0$ if $e\in\sigma$.
Hence $i^\sharp:F^\cdot(X)\goesto F^\cdot(Y)$ is surjective,
and $\ker(i^\sharp)=\span_\ZZ\{X_\sigma:\ e\in\sigma\in\Po(X)\}
=F^\cdot(X_e)(-1)$.  Fix an $\orn\in\Orn(X)$,
and denote by $\gamma$ its pullback along $i$.  Consider a generator
$R_\gamma(Y_\sigma,v)$ of $N_\gamma^\cdot(Y)$, so $\sigma\in\Po(Y)$ and
$v\in V(Y_\sigma)$.  Then
\showonq
i^\sharp(R_\orn (X_\sigma,v))
&=& \sum_{a\in X_\sigma(v)} s_\orn(a) i^\sharp(X_{\sigma\cup
\{\und{a}\}})\\
&=& \sum_{a\in Y_{\sigma}(v)} s_\gamma(a) Y_{\sigma\cup\{\und{a}\}}
=R_\gamma(Y_\sigma,v),
\showoffq
which shows that $i^\flat:N_\orn^\cdot(X)\goesto N_\gamma^\cdot(Y)$ is
surjective.  From the kernel-cokernel exact sequence, we see that
$i^\orn:K_\orn^\cdot(X)\goesto K_\gamma^\cdot(Y)$ is surjective and
$$0\lraw N_\orn^\cdot(X)\cap F^\cdot(X_e)(-1)\lraw
F^\cdot(X_e)(-1)\lraw\ker(i^\orn)\lraw 0$$
is exact.   If $e$ is a cut-edge then, as in Theorem $1.3$(a),
$F^\cdot(X_e)(-1)\subseteq N^\cdot_\orn(X)$, so $\ker(i^\orn)=0$;
since $\eta^\orn_{X,e}=0$ from Theorem $1.3$(a), this shows that
\showon
K^\cdot_\orn(X_e)(-1)\stackrel{\eta^\orn_{X,e}}{\lraw}
K^\cdot_\orn(X)\stackrel{i^\orn}{\lraw} K^\cdot_\gamma(Y)\lraw 0
\showoff
is exact in this case.  If $e$ is not a cut-edge then from Theorem
$1.3$(b), $\ker(\eta^\orn_{X,e})=0$, and so from (1.7) we find that
$N^\cdot_\orn(X)\cap F^\cdot(X_e)(-1)= N^\cdot_\orn(X_e)(-1)$. 
Injectivity of $\eta^\orn_{X,e}$ implies that (2.2) is exact in this
case as well.  The result follows since everything is compatible with
$\Psi(X)$ and the maps $\kappa_\gamma$ for $\gamma\in\Orn(Y)$.
\end{proof}

\begin{PROP} Let $h:Y\goesto X$ be a graph morphism such that
$h_V$ is bijective and $h_E$ is an elementary surjection.  Let $e$
be the unique edge of $X$ such that $\#h_E^{-1}(e)=2$.  Then
$$0\lraw K^\cdot(X)\stackrel{h^*}{\lraw} K^\cdot(Y)
\lraw K^\cdot(X_e)(-1)\oplus K^\cdot(X_e)(-2)\lraw 0$$
is exact. \end{PROP}
\begin{proof}
Identify $V(X)$ and $V(Y)$ via $h_V$, and consider $E(X)$ to be a
subset of $E(Y)$ by identifying each edge of $X$ with one member of
its preimage under $h_E$.  Let $c$ be the unique edge of $Y$ not in $X$,
and for $\sigma\in\Po(X)$ with $e\in\sigma$, let ${\sigma'}:=
\sigma\cup\{c\}\drop\{e\}$.
For all $\sigma\in\Po(X)$, $h^\sharp(X_\sigma)
:=Y_\sigma$ if $e\not\in\sigma$, and $h^\sharp(X_\sigma):=Y_\sigma+
Y_{{\sigma'}}$ if $e\in\sigma$.  Hence $h^\sharp$ is injective and
\showonq
\coker(h^\sharp) &=&
\frac{\ZZ\{Y_\sigma,Y_{{\sigma'}}:\ e\in\sigma\in\Po(X)\}}
{\ZZ\{Y_\sigma+Y_{{\sigma'}}:\ e\in\sigma\in\Po(X)\}}
\oplus\ZZ\{Y_\rho: \{c,e\}\subseteq\rho\in\Po(Y)\}\\
&\simeq& F^\cdot(X_e)(-1)\oplus F^\cdot(X_e)(-2).
\showoffq
Fix an $\orn\in\Orn(X)$, and denote by $\gamma$ its pullback along $h$.
Our next task is to determine the cokernel
of $h^\flat:N_\orn^\cdot(X)\goesto N_\gamma^\cdot(Y)$.  Consider any
$\sigma \in\Po(X)$ and $v\in V(X_\sigma)$; there are two cases.
If $e\not\in\sigma$ then $h^\flat(R_\orn(X_\sigma,v))=R_\gamma(Y_\sigma,v)$.
If $e\in\sigma$ then $h^\flat(R_\orn(X_\sigma,v))= R_\gamma(Y_\sigma,v)+
R_\gamma(Y_{{\sigma'}},v)$.  (In each case, the details of the
calculation differ slightly depending on whether $v$ is or is not
incident with $c$ in $Y_\sigma$.)  Therefore
\showonq
\coker(h^\flat) &=&
\frac{\span_\ZZ\{R_\gamma(Y_\sigma,v), :\ \sigma\in\Po(Y)\ 
\mathrm{and}\ \{e,c\}\cap\sigma\not\none\}}
{\span_\ZZ\{R_\gamma(Y_\sigma,v)+R_\gamma(Y_{{\sigma'}},v):\
e\in\sigma\in\Po(Y)\ \mathrm{and}\ c\not\in\sigma\}}\\
&\simeq& N_\orn^\cdot(X_e)(-1)\oplus N_\orn^\cdot(X_e)(-2),
\showoffq
an isomorphism compatible with the one above for $\coker(h^\sharp)$.
From the kernel-cokernel exact sequence we see that
$$0\lraw \ker(h^\orn)\lraw \coker(h^\flat)\stackrel{\phi}{\lraw}
\coker(h^\sharp)\lraw\coker(h^\orn)\lraw 0$$
is exact.  One checks that $\phi$ is injective, and it
follows that
$$0\lraw K_\orn^\cdot(X)\stackrel{h^\orn}{\lraw}K_\gamma^\cdot(Y)\lraw
K_\orn^\cdot(X_e)(-1)\oplus K_\orn^\cdot(X_e)(-2)\lraw 0$$
is exact.  The result follows for the usual reasons of compatibility.
\end{proof}

\begin{PROP} Let $g:Y\goesto X$ be a graph morphism such that
$g_V$ is an elementary surjection and $g_A$ is bijective.
Let $v$ be the unique vertex of $X$ such that $\#g_V^{-1}(v)=2$,
and let $g_V^{-1}(v)=\{u,w\}$.  Then
$$0\lraw \ker(g^*)\lraw K^\cdot(X)\stackrel{g^*}{\lraw}
K^\cdot(Y)\lraw 0$$
is exact, in which $\ker(g^*)\simeq J^\cdot(g):=\span_\ZZ
\{R_\gamma(Y_\sigma,u):\ \sigma\in \EuScript{Q}(g)\}$
and $\EuScript{Q}(g)$ is the set of all $\sigma\in\Po(Y)$ such that
the images of $u$ and $w$ in $Y_\sigma$ do not coincide.
\end{PROP}
\begin{proof}
Identify $V(X)\drop\{v\}$ with $V(Y)\drop\{u,w\}$ via $g_V$.
Identifying $E(X)$ with $E(Y)$ via $g_E$, we see that $g^\sharp:
F^\cdot(X)\goesto F^\cdot(Y)$ is the isomorphism defined by
$g^\sharp(X_\sigma):= Y_\sigma$ for all $\sigma\in\Po(X)$.  Fix an
$\orn\in\Orn(X)$ and let $\gamma$ be the pullback of $\orn$
along $g$.  To describe the cokernel of $g^\flat:N_\orn^\cdot(X)\goesto
N_\gamma^\cdot(Y)$, consider any $\sigma\in\Po(X)$ and $x\in V(X_\sigma)$.
If $x$ is not the image of $v$ in $V(X_\sigma)$ then
$g^\flat(R_\orn(X_\sigma,x))=R_\gamma(Y_\sigma,x)$.  If $x$ is the
image of $v$ in $X_\sigma$ and the images of $u$ and $w$ in $\sigma$
coincide, then $g^\flat(R_\orn(X_\sigma,x))=R_\gamma(Y_\sigma,u)$.
In the remaining case, $g^\flat(R_\orn(X_\sigma,x))=
R_\gamma(Y_\sigma,u)+R_\gamma(Y_\sigma,w)$.  Therefore,
$$\coker(g^\flat) =
\frac{\span_\ZZ\{R_\gamma(Y_\sigma,u),
R_\gamma(Y_\sigma,w):\ \sigma\in\EuScript{Q}(g)\}}
{\span_\ZZ\{R_\gamma(Y_\sigma,u)+R_\gamma(Y_\sigma,w):\ \sigma\in
\EuScript{Q}(g)\}}\simeq J^\cdot(g).$$
From the kernel-cokernel exact sequence we deduce that
$$0\lraw \ker(g^\orn)\lraw K_\orn^\cdot(X)\stackrel{g^\orn}{\lraw}
K_\gamma^\cdot(Y)\lraw 0$$
is exact, in which $\ker(g^\orn)\simeq\coker(g^\flat)$.
The proof is completed by compatibility of the maps, as usual.
\end{proof}

\section{Structure of Kirchhoff groups.}

For a graph $X$ and nonnegative integer $j$, let $d_j(X):=\rank\: K^j(X)$,
and collect these numbers as the coefficients of the
\emph{Poincar\'e polynomial} of $X$: $D_X(t):=d_0(X)+d_1(X)t+
\cdots+d_m(X)t^m$.  
We see in Theorem 3.1(a) that Kirchhoff groups are torsion-free, and so
the Poincar\'e polynomial of $X$ completely determines the structure of
$K^\cdot(X)$ as a graded abelian group.  However, the functoriality of
Kirchhoff groups induces a coalgebra structure on $K^\cdot(X)$, which
allows for the definition of more subtle algebraic invariants of graphs.

\begin{THM}  Let $X$ be a graph.\\
\textup{(a)}\ The Kirchhoff group $K^\cdot(X)$ is torsion-free.\\
\textup{(b)}\ Let $i:Y\goesto X$ be as in Proposition $2.4$, and
such that $e$ is not a cut-edge.  Then
$$D_X(t)= D_Y(t)+t\:D_{X_e}(t).$$
\textup{(c)}\ Let $h:Y\goesto X$ be as in Proposition $2.5$.  Then
$$D_Y(t)= D_X(t)+(t+t^2)D_{X_e}(t).$$
\textup{(d)}\ If $T_X(x,y)$ denotes the Tutte polynomial of $X$ then
$$D_X(t) = t^{n-k} T_X\left(\frac{1}{t},1+t\right),$$
in which $X$ has $n$ vertices and $k$ connected components.
\end{THM}
\begin{proof}
For part (a), we show that $K^\cdot(X)$ is torsion-free by induction
on $\#E$.  For the basis,
notice that $K^\cdot(X)=\ZZ$ for every graph $X$ with no edges.  For
a connected graph $\cutedge$ with two vertices and one
cut-edge, Proposition $2.4$ and Theorem $1.3$(a) imply that
$K^\cdot(\cutedge)= K^\cdot(\bullet\ \bullet)=\ZZ$.
Theorem 1.2 now implies that if every edge of $X$ is a cut-edge
(\emph{i.e.} $X$ is a forest) then $K^\cdot(X)=\ZZ$.  For the induction
step we may assume that $e\in E$ is not a cut-edge.  Theorem $1.3$(b)
and Proposition 2.4 show that
\showon
0 \lraw K^\cdot(X_e)(-1)\stackrel{\eta_{X,e}}{\lraw}
K^\cdot(X)\lraw K^\cdot(X\drop\{e\})\lraw 0
\showoff
is exact.  The homomorphism $\pi$ from Theorem $1.3$(b) shows that
this sequence is split, so that $K^\cdot(X)\simeq K^\cdot(X_e)(-1)\oplus
K^\cdot(X\drop\{e\})$.  By induction, both $K^\cdot(X_e)(-1)$ and
$K^\cdot(X\drop\{e\})$ are torsion-free, and so $K^\cdot(X)$ is torsion-free.

Parts (b) and (c) follow directly from $(3.1)$ and Proposition $2.5$,
respectively.

For part (d), first notice that $D_{\loopgph}(t)=1+t$ and
$D_{\cutedge}(t)=1$.  The function $\widetilde{D}_X(t):=t^{k-n}D_X(t)$ has
the property that if $e\in E$ is neither a loop nor a cut-edge, then
$$\widetilde{D}_X(t)=\widetilde{D}_{X\drop\{e\}}(t)+
\widetilde{D}_{X_e}(t),$$
as follows from part (b).  Since $T_X(x,y)$ is the universal
Tutte-Grothendieck invariant of the category of graphs (see Tutte \cite{Tu}
or Theorem 6.2.2 of Brylawski and Oxley \cite{BrOx}), it follows that
$\widetilde{D}_X(t)=T_X(\widetilde{D}_{\cutedge}(t),
\widetilde{D}_{\loopgph}(t))$.
The result now follows easily.
\end{proof}

\begin{CORO}  If $X$ and $Y$ are graphs with isomorphic graphic matroids,
then $K^\cdot(X)$ and $K^\cdot(Y)$ are isomorphic (but not naturally so).
\end{CORO}
\begin{proof}[First proof.]
If $X$ and $Y$ have isomorphic graphic matroids then $T_X(x,y)=T_Y(x,y)$
(since the Tutte polynomial is really a matroid invariant; see \cite{BrOx}).
If $L$ and $M$ are finitely generated torsion-free abelian
groups of the same rank then $L\simeq M$.
\end{proof}
\begin{proof}[Sketch of second proof.]
Whitney \cite{Wh} (see also Theorem $6.3.1$ of Oxley \cite{Ox}) proves that
two graphs have isomorphic graphic matroids if and only if they are
``$2$-isomorphic'' (see the above references for the definition).
If $X'$ is obtained from $X$ by splitting at a cut-vertex or by
merging two components to produce a cut-vertex then Theorem $1.2$ shows
that $K^\cdot(X')=K^\cdot(X)$. If $X'$ is obtained from $X$ by twisting
one component relative to a $2$-vertex-cut then fix an $\orn\in \Orn(X)$
and let $\gamma\in\Orn(X')$ be obtained from $\orn$ by changing the
orientation on all edges of the twisted component.  It is not difficult
to verify that $N^\cdot_\gamma(X')=N^\cdot_\orn(X)$, and so
$K^\cdot_\orn(X)=K^\cdot_\gamma(Y)$.  These equalities are compatible
with $\Psi(X)$ and $\Psi(X')$, and so $K^\cdot(X')= K^\cdot(X)$.
Now $X$ and $Y$ are $2$-isomorphic if and only if there is a graph
$X''$, obtained from $X$ by a finite sequence of operations $X\mapsto X'$
as above, such that $X''\simeq Y$.  The above remarks show that $K^\cdot(X'')
\simeq K^\cdot(X)$, but the isomorphism may depend nontrivially on a choice
of a sequence of operations $X\mapsto X'$ and thus is not natural. Finally, an
isomorphism $f:Y\goesto X''$ induces an isomorphism $f^*:K^\cdot(X'')\goesto
K^\cdot(Y)$; since this depends upon a choice for $f$ it is also not natural.
\end{proof}

For a graph $X=(V,A,o,t,\rev{\cdot})$ and a positive integer $r$,
define a graph $X^{(r)}:=(V',A',o',t',\rev{\cdot})$ as follows: it has
vertices $V':=V\times\{1,2,...,r\}$, arcs $A':=A\times\{1,2,...,r\}$, 
origin $o'(a,i):=(o(a),i)$, terminus $t'(a,i):=(t(a),i)$, and reversal
$\rev{(a,i)}:=(\rev{a},i)$.  This is a disjoint union of $r$ copies
of $X$, and is equipped with a natural graph morphism
$p_{(r)}:X^{(r)}\goesto X$ defined by projection onto the first coordinate.
By Theorem $1.2$, $K^\cdot(X^{(2)})=K^\cdot(X)\otimes K^\cdot(X).$
By functoriality of $K^\cdot$, there is a natural homomorphism
$\Delta:=p^*_{(2)}$
\showon
\Delta:K^\cdot(X)\lraw K^\cdot(X)\otimes K^\cdot(X).
\showoff
The automorphism of $X^{(2)}$ which exchanges the indices in the
second coordinate shows that $\Delta$ is cocommutative.
The two natural factorizations of $p_{(3)}$ through $p_{(2)}$ show that
$\Delta$ is coassociative.  With the projection $K^\cdot(X)\goesto
K^0(X)=\ZZ$ for counit, this gives $K^\cdot(X)$ the structure of a 
coalgebra over $\ZZ$.

\begin{PROP}  Let $f:Y\goesto X$ be a graph morphism.  Then
$f^*:K^\cdot(X)\goesto K^\cdot(Y)$ is coalgebra homomorphism.
\end{PROP}
\begin{proof}
Define $f^{(2)}:Y^{(2)}\goesto X^{(2)}$ by
$f^{(2)}_V:= f_V\times\mathrm{id}_{\{1,2\}}$ and
$f^{(2)}_A:= f_A\times\mathrm{id}_{\{1,2\}}$.  The diagram of
graph morphisms
$$\begin{array}{ccc}
Y^{(2)} & \stackrel{f^{(2)}}{\longrightarrow} & X^{(2)}\\
p_{(2)}\downarrow\ \ & &\ \ \downarrow p_{(2)}\\
Y & \stackrel{f}{\longrightarrow} & X\end{array}$$
commutes.  Functoriality of $K^\cdot$ now implies the result,
since $(f^{(2)})^*=f^*\otimes f^*$.
\end{proof}

We remark that the abelian group isomorphisms constructed in the second
proof of Corollary $3.2$ are in fact coalgebra isomorphisms.

\setcounter{section}{3}
\section{Circulation algebras of graphs.}

Next, we dualize the coalgebra structure on $K^\cdot(X)$ to 
obtain functorial algebra structures on $\Hom(K^\cdot(X),B)$ for any
ring $B$, give an explicit description of these algebras when $B$
is commutative, and use this to derive the results on the ranks
$d_j(X)$ presented in the Introduction.

For a ring $B$ with unit $1$, let $\Flow_\cdot(X,B):=\Hom(K^\cdot(X),(B,+))$.
The ring structure of $B$ is a map $B\otimes B\goesto B$,
which induces a natural homomorphism
\showonq
\Flow_\cdot(X,B)\otimes\Flow_\cdot(X,B) &
\stackrel{\sim}{\lraw} &
\Hom(K^\cdot(X)\otimes K^\cdot(X),B\otimes B)\\
& \lraw & \Hom(K^\cdot(X)\otimes K^\cdot(X),B).
\showoffq
(The first map is an isomorphism.)  This factors through
$\Flow_\cdot(X,B)\otimes_B\Flow_\cdot(X,B)$, giving
$$\Flow_\cdot(X,B)\otimes_B\Flow_\cdot(X,B)\lraw
\Hom(K^\cdot(X)\otimes K^\cdot(X),B).$$
Dualizing (3.2), we obtain 
$$\Hom(K^\cdot(X)\otimes K^\cdot(X),B)\lraw \Flow(X,B).$$
The composition of these gives a natural homomorphism 
$$\Flow_\cdot(X,B)\otimes_B\Flow_\cdot(X,B) \lraw\Flow_\cdot(X,B),$$
and this defines a $B$-algebra structure on $\Flow_\cdot(X,B)$.
Since $\Delta$ is cocommutative, $\Flow_\cdot(X,B)$ is commutative
if and only if $B$ is commutative.  From the functoriality of $\Hom$
and $K^\cdot$ and Proposition $3.3$,
it follows that $\Flow_\cdot$ is a functor which is covariant in both
of its arguments.  In particular, a graph morphism $f:Y\goesto X$
induces a $B$-algebra homomorphism $f_*:\Flow_\cdot(Y,B)\goesto
\Flow_\cdot(X,B)$.

The \emph{$B$-circulation algebra} of a graph $X$ is
$\Flow_\cdot(X,B)$;  it is a graded $B$-algebra which is a
finitely generated $B$-module.  As mentioned above,
for any $\varphi\in\Flow_\cdot(X,B)$ the $j$-th graded component
of $\varphi$ is a coherent family of $B$-valued flows on
$\{X_\sigma:\ \sigma\in\Po_{j-1}(X)\}$.  An element of
$\Flow_\cdot(X,B)$ will be called a \emph{$B$-circulation on $X$};
the \emph{$B$-flows in $X$} are precisely the $B$-circulations on $X$ which
are in $\Flow_1(X,B)$.
Each $B$-circulation on $X$ is a homomorphism $\varphi:
F^\cdot_\Omega(X)\goesto(B,+)$ which annihilates
$N^\cdot_\Omega(X)$.  Thus $\Flow_\cdot(X,B)$ is naturally a $B$-submodule
of $\Hom(F^\cdot_\Orn(X),B)$.  Fix an $\orn\in\Orn$.  For each
$\varphi\in\Flow_\cdot(X,B)$, the composition
$$F^\cdot(X)\lraw F^\cdot(X)/N_\orn^\cdot(X)\stackrel{\kappa_\orn}{\lraw}
K^\cdot(X)\stackrel{\varphi}{\lraw}B$$
is an element of $\Hom(F^\cdot(X),B)$.  This defines a monomorphism
$$\lambda_\orn:\Flow_\cdot(X,B)\goesto\Hom(F^\cdot(X),B)$$
which identifies $\Flow_\cdot(X,B)$ with the submodule of
$\Hom(F^\cdot(X),B)$ annihilating $N^\cdot_\orn(X)$;  it will be
convenient to refer to this identification as \emph{coordinatization of
$\Flow_\cdot(X,B)$ by $\orn$}.  The basis of $\Hom(F^\cdot(X),B)$
dual to $\{X_\sigma:\ \sigma\in\Po(X)\}$ will be denoted by
$\{\xi_\sigma:\ \sigma\in\Po(X)\}$.

There is a
comultiplication $\hat{\Delta}$ defined on $F^\cdot(X)$ by 
$$\hat{\Delta}(X_\sigma):= \sum_{\tau\subseteq\sigma} X_\tau\otimes
X_{\sigma\drop\tau}$$
for each $\sigma\in\Po(X)$.  Dualizing gives a multiplication on
$\Hom(F^\cdot(X),B)$:  for any $\varphi,\theta\in\Hom(F^\cdot(X),B)$
and any $\sigma\in\Po(X)$,
\showon
(\varphi\cdot\theta)(X_\sigma)= \sum_{\tau\subseteq\sigma}
\varphi(X_\tau)\theta(X_{\sigma\drop\tau}).
\showoff
A graph morphism $f:Y\goesto X$ induces a $B$-algebra homomorphism
$f_\sharp:\Hom(F^\cdot(Y),B)\goesto\Hom(F^\cdot(X),B)$ by
$f_\sharp(\varphi):=\varphi\circ f^\sharp$, and this is functorial.
Restriction of the multiplication on $\Hom(F^\cdot(X),B)$ to the
annihilator of $N_\orn^\cdot(X)$ coincides with the multiplication on
$\Flow_\cdot(X,B)$ via the map $\lambda_\orn$.

The \emph{positive ideal} of $\Hom(F^\cdot(X),B)$ is $\Hom(F^+(X),B)
=$ $\bigoplus_{j>0}\Hom(F^j(X),B)$, and for circulation algebras
$\Flow_+(X,B):=\bigoplus_{j>0}\Flow_j(X,B)$.  For each
$\varphi\in\Hom(F^+(X),B)$ there is a nonnegative integer $r$ such
that $\varphi^r=0$;  the
\emph{nilpotence} np$(\varphi)$ of such a $\varphi$ is the greatest
integer $n$ such that $\varphi^n\neq 0$.
For $\varphi\in\Hom(F^\cdot(X),B)$, define the \emph{support of $\varphi$}
to be $\supp(\varphi):=\{\sigma\in\Po(X):\ \varphi(X_\sigma)\neq 0\}$.
For $\orn,\eps\in\Orn$, $\varphi\in\Flow_\cdot(X,B)$, and
$\sigma\in\Po(X)$, notice that $\lambda_\eps(\varphi)(X_\sigma)=
(-1)^{m_\sigma(\eps,\orn)} \lambda_\orn(\varphi)(X_\sigma)$;  that is,
$\lambda_\eps(\varphi)= \lambda_\orn(\varphi)\circ\psi^\eps_\orn$.
Thus, the set $\supp(\lambda_\orn(\varphi))$ does not depend on the
choice of $\orn\in\Orn$, so each $B$-circulation on $X$ also has a
well-defined support.
\begin{LMA}  Let $X$ be a graph, let $B$ be a commutative ring, and let
$\varphi\in\Hom(F^+(X),B)$.\\
\textup{(a)}\ For every $\sigma\in\Po(X)$ and nonnegative integer $r$,
$$\varphi^r(X_\sigma)=r!\sum_{\{\tau_1,...,\tau_r\}}^\sigma
\prod_{i=1}^r\varphi(X_{\tau_i}),$$
in which the sum is over all partitions of the set $\sigma$ into
exactly $r$ pairwise disjoint nonempty blocks.\\
\textup{(b)}\ If $\varphi\in\Hom(F^1(X),B)$ then \textup{np}$(\varphi)\leq
\#\supp(\varphi)$.\\
\textup{(c)}\ If $\varphi\in\Hom(F^1(X),B)$ and $B$ is an integral domain
with $\cha(B)=0$, then \textup{np}$(\varphi)=\#\supp(\varphi)$.\\
\textup{(d)}\ If $\varphi\in\Hom(F^1(X),B)$ and $B$ is an integral domain
with $\cha(B)=p>0$, then \textup{np}$(\varphi)=\min\{p-1,\#\supp(\varphi)\}$.
\end{LMA}
\begin{proof} From $(4.1)$ we obtain
$$\varphi^r(X_\sigma)=\sum_{(\tau_1,...,\tau_r)}^\sigma\varphi(X_{\tau_1})
\varphi(X_{\tau_2})\cdots\varphi(X_{\tau_r}),$$
in which the sum is over all ordered partitions of the set $\sigma$
into exactly $r$ pairwise disjoint sets.  If any of the $\tau_i=\none$
then the corresponding term is zero, since $\varphi(X_\none)=0$.
Since $B$ is commutative, the $r!$ terms corresponding to the same
unordered partition all have the same value, proving part (a).

If $\varphi\in\Hom(F^1(X),B)$ then $\varphi(X_\tau)=0$
unless $\tau\in\Po_1(X)$.  It follows that for all $\sigma\in\Po_j(X)$,
$$\varphi^r(X_\sigma)=r!\delta_{jr}\prod_{e\in\sigma}\varphi(X_e),$$
in which $\delta_{jr}$ is the Kronecker delta function.  If $j>
\#\supp(\varphi)$ then every $\sigma\in\Po_j(X)$ contains some $e\not\in
\supp(\varphi)$, and it follows that $\varphi^j=0$.  This proves (b).
If $B$ is an integral domain and $\supp(\varphi)=\{e_1,...,e_r\}$
then for all $0\leq j\leq r$, $\varphi^j(X_{\{e_1,...,e_j\}})\neq 0$
if and only if $j!\neq 0$ in $B$.  Since $\cha(B)$ is either zero or
prime, parts (c) and (d) follow.
\end{proof} 

When $B$ is commutative, the \emph{exponential function}
$$\exp:\Hom(F^+(X),B)\goesto\Hom(F^\cdot(X),B)$$
is defined for each $\varphi\in\Hom(F^+(X),B)$ by putting
\showon
\exp(\varphi)(X_\sigma):=\sum_{r=0}^\infty\sum_{\{\tau_1,...,\tau_r\}}^\sigma
\prod_{i=1}^r\varphi(X_{\tau_i})
\showoff
for all $\sigma\in\Po(X)$, with notation as in Lemma $4.1$(a).
If $B$ is a $\QQ$-algebra then one may also use the formal power series
$\exp(\varphi)=\sum_{r=0}^\infty\varphi^r/r!$
because of Lemma $4.1$(a).  For $\varphi\in\Hom(F^1(X),B)$ we use the
notation
$$\exp(\varphi)=\varphi^{<0>}+\varphi^{<1>}+\varphi^{<2>}+\cdots
+\varphi^{<j>}+\cdots$$
in which $\varphi^{<j>}$ is the $j$-th graded component of $\exp(\varphi)$,
so $\varphi^{<0>}=1$ and $\varphi^{<1>}=\varphi$.
\begin{LMA}  Let $X$ be a graph and let $B$ be a commutative ring.\\
\textup{(a)}\  For $\varphi,\theta\in\Hom(F^+(X),B)$, we have
$\exp(\varphi+\theta)=\exp(\varphi)\cdot\exp(\theta)$.\\
\textup{(b)}\ If $f:Y\goesto X$ is a graph morphism then 
$\exp\circ f_\sharp=f_\sharp\circ\exp$.\\
\textup{(c)}\ For any $\orn\in\Orn$, we have
$\exp(\lambda_\orn(\Flow_+(X,B)))\subseteq\lambda_\orn(\Flow_\cdot(X,B)).$
\end{LMA}
\begin{proof}
Part (a) is a straightforward calculation from $(4.1)$ and $(4.2)$.
For (b), let $\varphi\in\Hom(F^\cdot(Y),B)$ and $\sigma\in\Po(X)$, and
calculate that
\showonq
\exp(f_\sharp(\varphi))(X_\sigma) &=& \exp(\varphi\circ f^\sharp)(X_\sigma)\\
&=& \sum_{r=0}^\infty \sum_{\{\tau_1,...,\tau_r\}}^\sigma\prod_{i=1}^r
\sum_{\rho_i\in T_f(\tau_i)}\varphi(Y_{\rho_i})\\
&=& \sum_{\rho\in T_f(\sigma)}\sum_{r=0}^\infty
\sum_{\{\tau_1,...,\tau_r\}}^\rho\prod_{i=1}^r\varphi(Y_{\tau_i})\\
&=& \sum_{\rho\in T_f(\sigma)}\exp(\varphi)(Y_\rho)\\
&=& \exp(\varphi)(f^\sharp(X_\sigma))=f_\sharp(\exp(\varphi))(X_\sigma).
\showoffq
For part (c) let $\varphi\in\lambda_\orn(\Flow_+(X,B))$, consider
a relation $R_\orn(X_\sigma,v)$ in $N^\cdot_\orn(X)$, and calculate that
\showonq
& & \exp(\varphi)(R_\orn(X_\sigma,v))\\
&=& \exp(\varphi)\left(\sum_{a\in X_\sigma(v)} s_\orn(a)
    X_{\sigma\cup\{\und{a}\}}\right)
= \sum_{a\in X_\sigma(v)} s_\orn(a) \sum_{r=0}^\infty \sum_{\{\tau_1,
...,\tau_r\}}^{\sigma\cup\{\und{a}\}}\prod_{i=1}^r\varphi(X_{\tau_i})\\
&=& \sum_{r=0}^\infty\sum_{\{\tau_1,...,\tau_r\}}^\sigma\left[\sum_{h=1}^r
\left(\prod_{i\neq h}\varphi(X_{\tau_i})\right)\sum_{a\in
X_\sigma(v)} s_\orn(a) \varphi(X_{\tau_h\cup\{\und{a}\}})\right.\\
& & + \left.
\left(\prod_{i=1}^r\varphi(X_{\tau_i})\right)\sum_{a\in
X_\sigma(v)} s_\orn(a) \varphi(X_{\und{a}})\right]\\
&=& \sum_{r=0}^\infty \sum_{\{\tau_1,...,\tau_r\}}^\sigma\left[
\sum_{h=1}^r\left(\prod_{i\neq h}\varphi(X_{\tau_i})\right)
\varphi\left(\sum_{u\in U(h)} R_\orn(X_{\tau_h},u)\right)\right.\\
& & + \left.
\left(\prod_{i=1}^r\varphi(X_{\tau_i})\right)
\varphi\left(\sum_{u\in U} R_\orn(X,u)\right)\right]= 0.
\showoffq
In the last expression, for each partition $\{\tau_1,...,\tau_r\}$ of
$\sigma$, $U(h)$ denotes the set of vertices of $X_{\tau_h}$ which are
contracted to the vertex $v$ of $X_\sigma$ (for each $1\leq h\leq r$),
and $U$ denotes the set of vertices of $X$ which are contracted to the
vertex $v$ of $\sigma$. (The penultimate equality follows from Lemma $1.1$.)
This proves part (c).
\end{proof}

For a maximal forest $Y$ of $X$ and an edge $e$ of $X$ not in $Y$,
the subgraph $Y\cup\{e\}$ contains a unique cycle $C$, the \emph{
fundamental cycle of $e$ with respect to $Y$}.  Let $a\in e$
be one of the arcs in $e$, let $\orn\in\Orn$ be such that $s_\orn(a)=1$,
and coordinatize $\Flow_\cdot(X,B)$ by $\orn$.  There is a unique flow
$\beta\in\Flow_1(C,B)\subseteq\Flow_1(X,B)$ such that $\beta(X_e)=1$,
and it is independent of the choice of $\orn\in\Orn$ (subject to
$s_\orn(a)=1$).  This is the \emph{basic flow of $a$ with respect
to $Y$}.

It is convenient to standardize some notation which will be
used repeatedly in what follows.\\
\textbf{(*)}\ \
Coordinatize $\Flow_\cdot(X,B)$ by some $\orn\in\Orn$, and fix a maximal
forest $Y$ contained in $X$.  Let $\EE:=E(X)\drop E(Y)$, and for each
$c\in\EE$ let $\beta_c$ be the basic flow of $\orn(c)$ with respect to $Y$.
For $\jj:\mathcal{E}\goesto\NN$ define $\BETA^{<\jj>}:=
\prod_{c\in\mathcal{E}}\beta_c^{<j(c)>}$.\\
\textbf{(**)}\ \ 
Continuing from (*), fix an $e\in\EE$ which is not a loop.
Then the image of $Y$ in $X_e$ contains a unique cycle $C_e$, and this
is the image in $X_e$ of the unique cycle $C$ in $Y\cup\{e\}$.  We
may choose $\orn\in\Orn$ so that for every vertex $v\in V(C)$ there
is a unique edge $c\in E(C)$ such that $v=t(\orn(e))$.  (Changing the
orientation on a subset of $E(Y)$ does not affect any of the basic
flows $\{\beta_c:\ c\in\EE\}$.)  Let $\zeta:=E(C_e)$, and let $e'$ be
any edge in $\zeta$.  Then $Y':=Y\drop\{e'\}$
is a maximal forest in $X_e$. Let $\EE':=E(X_e)\drop E(Y')$ and for each
$c\in\EE'$ let $\alpha_c$ be the basic flow of $\orn(c)$ with respect
to $Y'$.  Let $\pi:K^\cdot(X)\goesto K^\cdot(X_e)(-1)$ be the map
constructed from $C$ in the proof of Theorem $1.3$(b), and let
$\pi^*:\Flow_\cdot(X_e,B)(-1) \goesto\Flow_\cdot(X,B)$ be the map
induced from $\pi$ by duality.

We require one more lemma for the proofs of the main results of this
section.  For $e\in E$ define $p_e:\Hom(F^\cdot(X),B)\goesto
\Hom(F^\cdot(X_e),B)$ by putting $p_e(\theta)((X_e)_\sigma)
:= \theta(X_\sigma)$ for all $\theta\in\Hom(F^\cdot(X),B)$ and
$\sigma\in\Po(X_e)$.

\begin{LMA} Let $X$ be a graph, let $e\in E$, and let $B$ be a commutative
ring.\\
\textup{(a)}\ The map $p_e$ is a $B$-algebra homomorphism.\\
\textup{(b)}\ For any $\orn\in\Orn$, we have
$p_e(\lambda_\orn(\Flow_\cdot(X,B)))\subseteq
\lambda_\orn(\Flow_\cdot(X_e,B)).$\\
\textup{(c)}\  With the notation \textup{(**)}, $p_e(\beta_e)=\alpha_{e'}$
and if $c\in\EE\drop\{e\}$ then $p_e(\beta_c)=\alpha_c+s_c\alpha_{e'}$,
in which $s_c:=\beta_c(X_{e'})\in\{-1,0,1\}$.\\
\textup{(d)}\  With the notation \textup{(**)}, $(\pi^*\circ p_e)
(\BETA^{<\jj>})=\BETA^{<\jj+\delta_e>}$ for all $\jj:\EE\goesto\NN$,
in which $\delta_e:\EE\goesto\NN$ is defined by
$\delta_e(c):=\tv[c=e]$ for all $c\in\EE$.
\end{LMA}
\begin{proof}
Parts (a) and (b) are straightforward.  For part (c),  if $c\in\EE$ then
let $C(c)$ be the fundamental cycle of $c$ with respect to $Y$, and if
$c\in\EE'$ then let $C'(c)$ be the fundamental cycle of $c$ with respect
to $Y'$.  Then $C'(e')=C_e$ as defined in (**), and
for $c\in\EE'\drop\{e'\}$, if $e'\not\in E(C(c))$ then $C'(c)=C(c)$, and
if $e'\in E(C(c))$ then $C'(c)$ has edge-set the symmetric difference of
$E(C(c))$ and $E(C_e)$.  Also, for $c,g\in\EE$, $\beta_c(X_g)=0$ if $g$
is not an edge of $C(c)$, and otherwise $\beta_c(X_g)=\pm 1$
according to whether $\orn(c)$ and $\orn(g)$ are directed in the same
or in opposite directions around $C(c)$.  From this part (c) follows.
For part (d), notice for $\jj:\EE\goesto\NN$ and $\sigma\in\Po(X)$ that
$$\BETA^{<\jj>}(X_\sigma)=\sum_{(\tau(c):\ c\in\EE)}^\sigma
\prod_{c\in\EE}\prod_{g\in\tau(c)}\beta_c(X_g),$$
in which the sum is over the set of all $\EE$-indexed partitions
$(\tau(c):\ c\in\EE)$ of $\sigma$ into disjoint subsets
such that $\#\tau(c)=j(c)$ for all $c\in\EE$.
Now for all $\jj:\EE\goesto\NN$ and $\sigma\in\Po(X)$,
\showonq
(\pi^*\circ p_e)(\BETA^{<\jj>})(X_\sigma)
&=& p_e(\BETA^{<\jj>})(\pi(X_\sigma))\\
&=&\left\{\begin{array}{ll}
p_e(\BETA^{<\jj>})((X_e)_{\sigma\drop\{e\}}) & \mathrm{if}\ e\in\sigma,\\
\sum_{c\in\sigma\cap\zeta} p_e(\BETA^{<\jj>})((X_e)_{\sigma\drop\{c\}})
& \mathrm{if}\ e\not\in\sigma,\end{array}\right.\\
&=&\left\{\begin{array}{ll}
\BETA^{<\jj>}(X_{\sigma\drop\{e\}}) & \mathrm{if}\ e\in\sigma,\\
\sum_{c\in\sigma\cap\zeta}\BETA^{<\jj>}(X_{\sigma\drop\{c\}})
& \mathrm{if}\ e\not\in\sigma,
\end{array}\right.\\
&=&\left\{\begin{array}{ll}
\BETA^{<\jj+\delta_e>}(X_{\sigma}) & \mathrm{if}\ e\in\sigma,\\
\BETA^{<\jj+\delta_e>}(X_{\sigma}) & \mathrm{if}\ e\not\in\sigma.
\end{array}\right.
\showoffq
The last equality in case $e\not\in\sigma$ follows from (4.1) and
the fact that $\beta_e=\sum_{c\in\zeta}\xi_c$.
Hence $(\pi^*\circ p_e)(\BETA^{<\jj>})=\BETA^{<\jj+\delta_e>}$ for all
$\jj:\EE\goesto\NN$.
\end{proof}

\begin{THM}  Let $X$ be a graph, let $B$ be a commutative ring,
and adopt the notation \textup{(*)}.\\
\textup{(a)}\  Then $\{\BETA^{<\jj>}:\ \jj:\EE\goesto\NN\}$
spans $\Phi_\cdot(X,B)$ as a $B$-module.\\
\textup{(b)}\ If $e\in\EE$ is not a loop then $p_e$ is surjective.
\end{THM}
\begin{proof}  We prove (a) and (b) together by induction on $\#E(X)$.
As basis of induction we have the case that $X$ is a forest;  then
$Y=X$ and $\EE=\none$ and $\Flow_\cdot(X,B)=B$ is spanned by
$\{1\}$ as a $B$-module, conforming with the statement to be proved.
For the induction step, let $e\in \EE$, so $e$ is not a cut-edge of $X$.
If $e$ is a loop then the dual of Theorem $1.2$ implies that
\showon
\Flow_\cdot(X,B)=\Flow_\cdot(X\drop\{e\},B)\otimes_B\Flow_\cdot(\loopgph,B),
\showoff
and the induction argument is easy.
In the remaining case, adopt the notation (**).  Dualizing $(3.1)$
we see that
$$0\lraw\Flow_\cdot(X\drop\{e\},B)\stackrel{i_*}{\lraw}
\Flow_\cdot(X,B)\stackrel{\eta_{X,e}^*}{\lraw}\Flow_\cdot(X_e,B)(-1)
\lraw 0$$
is exact, and the map $\pi^*$ shows that
the sequence is split.  Thus, any $\varphi\in\Flow_\cdot(X,B)$ can be
written in the form $\varphi=i_*(\varphi'')+\pi^*(\varphi')$ with
$\varphi''\in\Flow_\cdot(X\drop\{e\},B)$ and $\varphi'\in
\Flow_\cdot(X_e,B)(-1)$.  By induction, $\Flow_\cdot(X\drop\{e\},B)$
is spanned by $\{\BETA^{<\jj>}:\ \jj:\EE\drop\{e\}\goesto\NN\}$.  Since
$i_*(\beta_c)=\beta_c$ for all $c\in \EE\drop\{e\}$ and $i_*$
is a $B$-algebra homomorphism, the term $i_*(\varphi'')$ is of the
required form.

To understand the term $\pi^*(\varphi')$, the induction hypothesis implies
that $\{\boldsymbol{\alpha}^{<\mathbf{h}>}:\ \mathbf{h}:\EE'\goesto\NN\}$
spans $\Flow_\cdot(X_e,B)$ as a $B$-module. From Lemma 4.3(c) it follows for
all $h\in\NN$ that ${\alpha_{e'}}^{<h>}=p_e(\beta_e^{<h>})$ and if
$c\in\EE'\drop\{e'\}$ then ${\alpha_c}^{<h>}=
p_e\left((\beta_c-s_c\beta_e)^{<h>}\right)=
p_e\left(\sum_{j=0}^h \beta_c^{<h-j>}\cdot (-s_c)^j\beta_e^{<j>}\right)$.
Since $\beta_e^{<i>}\cdot\beta_e^{<j>}=\binom{i+j}{j}
\beta_e^{<i+j>}$ it follows for all $\mathbf{h}:\EE'\goesto\NN$ that
$\ALPHA^{<\mathbf{h}>}$ is in the span over $B$ of the elements
$p_e(\BETA^{<\jj>})$ for $\jj:\EE\goesto\NN$.
This proves (b), and so there is a (finite) expression
$$\varphi'=\sum_{\jj:\EE\goesto\NN} b(\jj)p_e(\BETA^{<\jj>})$$
with coefficients $b(\jj)$ in $B$.  Therefore, from Lemma 4.3(d) we
conclude that
$$\pi^*(\varphi')=\sum_{\jj:\EE\goesto\NN} b(\jj)
\BETA^{<\jj+\delta_e>},$$
completing the induction step, and the proof of (a).
\end{proof}

For example, applying Theorem $4.4$ to the $n$-vertex cycle $\Cyc_n$
and $B=\ZZ$ we see that
$$\Flow_\cdot(\Cyc_n,\ZZ)=\span_\ZZ\{1,\beta,\beta^{<2>},...,\beta^{<n>}\},$$
in which $\beta$ is the basic flow with respect to some arc and maximal
tree in $\Cyc_n$ (since np$(\beta)=n$, by Lemma $4.1$(c)).  This is
isomorphic to the quotient of the \emph{divided power algebra} $\ZZ_![x]:=
\ZZ[x^j/j!:\ j\in\NN]$ (a non-Noetherian subring of $\QQ[x]$) by the
ideal $(x^j/j!:\ j>n)$.  The notation $x^{<j>}:=x^j/j!$ for $j\in\NN$
will be convenient.  Notice that $d/dx$ and $\int dx$ are endomorphisms
of $\ZZ_![x]$ which are surjective and injective, respectively.

\begin{CORO}  Let $X$ be a graph, let $B$ be a commutative 
$\QQ$-algebra, and adopt the notation \textup{(*)}. 
Then $\{\beta_c:\ c\in\EE\}$ generates $\Flow_\cdot(X,B)$ as
a $B$-algebra.\end{CORO}
\begin{proof}
This follows immediately from Theorem 4.4, since $\beta_c^{<j>}=
\beta_c^j/j!$ for all $c\in\EE$ and $j\in\NN$.
\end{proof}
\begin{CORO} Let $X$ be a graph with $n$ vertices, $m$ edges, $k$
connected components, and $\ell$ cut-edges, and let
$D_X(t)=d_0+d_1t+\cdots+d_mt^m$.  Then $d_0=1$, $d_1=m-n+k$,
$d_{m-\ell}=1$, $d_j\neq 0$ if and only if  $0\leq j\leq m-\ell$,
and with the notation of $(0.2)$, if $1\leq j\leq m-\ell-1$ then
$d_{j+1}\leq\psi_j(d_j)$.
\end{CORO}
\begin{proof}
Adopt the notation (*) with $B=\QQ$.
We have $d_j=\dim_\QQ\Flow_j(X,\QQ)$ for all $0\leq j\leq m$.
The equality $d_0=1$ is clear, since $\Flow_0(X,\QQ)=\QQ$.
From Theorem $4.4$, the basic flows $\{\beta_c:\ c\in\EE\}$ span
$\Flow_1(X,\QQ)$, and since $\beta_c(X_g)=\delta_c(g)$ for all
$c,g\in\EE$, these flows are linearly independent.  Hence
$d_1=\#\EE=m-n+k$, as is well-known.  Since $\Flow_\cdot(\cutedge,\QQ)=\QQ$,
Theorem $1.2$ (dualized) allows us to reduce to the case that $X$ has no
cut-edges, and implies in general that $d_j=0$ for $m-\ell<j\leq m$.
When $X$ has no cut-edges, $F^m(X)$ is spanned by $X_E$ and for any
$\sigma\in\Po_{m-1}(X)$ the graph $X_\sigma$ is isomorphic with $\loopgph$;
thus, for the unique $v\in V(X_\sigma)$ we have $R_\orn(X_\sigma,v)=0$.
Hence $N^m(X)=0$, and it follows that $K^m(X)\simeq\ZZ$, and so $d_m=1$
(in case $X$ has no cut-edges).  Since $\Flow_1(X,\QQ)$ generates
$\Flow_\cdot(X,\QQ)$ as a $\QQ$-algebra, $d_j\neq 0$ for all $0\leq
j\leq m-\ell$, and the inequalities $d_{j+1}\leq \psi_j(d_j)$ follow
from Macaulay's Theorem (see Theorems II.2.2 and II.2.3 of Stanley
\cite{St}).
\end{proof}
\begin{CORO} Let $X$ be a graph with $m$ edges and $\ell$ cut-edges,
and let $D_X(t):=d_0+d_1t+\cdots+d_{m-\ell}t^{m-\ell}$.  Adopt
the notation \textup{(*)} with $B=\QQ$, and let $r(c):=\#\supp(\beta_c)$
for all $c\in\EE$.  Then, for $0\leq j\leq m-\ell$,
$$d_j\leq [t^j]\prod_{c\in\EE}\left(1+t+t^2+\cdots+t^{r(c)}\right).$$
\end{CORO}
\begin{proof}
For $c\in\EE$, Lemma $4.1$(c) shows that $\beta_c^{1+r(c)}=0$, and
so $\{\BETA^\jj:\ 0\leq j(c)\leq r(c)\ 
\mathrm{for}\ \mathrm{all}\ c\in\EE\}$ spans $\Flow_\cdot(X,\QQ)$ over $\QQ$.
If these monomials are linearly independent then the Poincar\'e
polynomial on the right side of the inequality is attained; in any case
the dimensions $d_j$ are bounded above as described.
\end{proof}
In fact, a common generalization of Corollaries $4.6$ and $4.7$ can
be obtained by applying the Clements-Lindstr\"om Theorem \cite{CL},
but we will not state the resulting inequalities explicitly here.

Next, for a commutative ring $B$, we present the algebra $\Flow_\cdot(X,B)$
explicitly as a quotient of a (multivariate) divided power algebra.  Adopt
the notation (*), let
$$\ZZ_![\xx]:=\bigotimes_{c\in\EE}\ZZ_![x_e],$$
in which $\xx:=\{x_c:\ c\in\EE\}$ are indeterminates algebraically independent
over $\QQ$, and let $B_![\xx]:=B\otimes\ZZ_![\xx]$.
For $\jj:\EE\goesto\NN$ we use the notation
$\xx^{<\jj>}$ for $\prod_{c\in\EE}x_c^{<j(c)>}$; these elements form a basis
for $B_![\xx]$ as a free $B$-module.  For a homogeneous
linear polynomial $P(\xx)=\sum_{c\in\EE}b_cx_c\in B_![\xx]$ and any $r\in\NN$,
we define $P^{<r>}(\xx)$ by
$$P^{<r>}(\xx):=\sum_\jj \prod_{c\in\EE}b_c^{j(c)}x_c^{<j(c)>},$$
with the sum over all $\jj:\EE\goesto\NN$ such that $\sum_{c\in\EE}j(c)=r$;
when $B$ is a $\QQ$-algebra this is consistent with the multinomial 
theorem expansion of $P(x)^{r}/r!$.
Define a $B$-linear homomorphism $\phi:B_![\xx]\goesto\Phi_\cdot(X,B)$
by putting $\phi(\xx^{<\jj>}):=\BETA^{<\jj>}$ for all $\jj:\EE\goesto\NN$;
by Theorem 4.4, $\phi$ is surjective.  It is not difficult to verify 
that $\phi$ is in fact a ring homomorphism.  Consider any
$\theta\in\Flow_1(X,B)$, and let $r\geq\#\supp(\theta)$. The linear polynomial
$$P_\theta(\xx):=\sum_{c\in\EE}\theta(X_c)x_c$$
is such that $P_\theta^{<1+r>}(\xx)$ is in the kernel of $\phi$, as follows from
$(4.2)$ and the fact that $\theta=\sum_{c\in\EE}\theta(X_c)\beta_c$.
Let $\mathcal{C}$ be the subset of $\Flow_1(X,B)$ consisting of those
flows $\theta$ such that $\supp(\theta)$ is a cycle and for each $c\in E$,
$\theta(X_c)\in\{-1,0,1\}$.  Define an ideal in $B_![\xx]$ by
$$I(X,Y,\orn):=(P_\theta^{<1+r>}(\xx):\
\theta\in\mathcal{C}\ \mathrm{and}\ r\geq\#\supp(\theta)).$$

\begin{THM} Let $X$ be a graph, let $B$ be a commutative ring,
and adopt the notation \textup{(*)}.  Then 
$$\Flow_\cdot(X,B)\simeq\frac{B_![\xx]}{I(X,Y,\orn)}.$$
\end{THM}
\begin{proof}
We proceed by induction on $\#E(X)$, using as basis the case that 
$X$ is a forest, in which case the claim is trivial.  For the
induction step, let $e\in\EE$.  If $e$ is a loop of $X$ then $(4.3)$
holds and the induction step is easy.  Otherwise, adopt the notation
(**).  Since $e$ is not a cut-edge of $X$, we have a commutative diagram
$$\begin{array}{ccccccccc}
0& \lraw&  B_![x_c:\ c\in\EE\drop\{e\}]& \stackrel{\tilde{\iota}}{\lraw}&
B_![\xx]& \stackrel{\tilde{\eta}}{\lraw}& B_![\yy](-1)&\lraw& 0\\
& & \downarrow\phi''& & \downarrow\phi& & \downarrow\phi' & & \\
0& \lraw&  \Flow_\cdot(X\drop\{e\},B)& \stackrel{i_*}{\lraw}& 
\Flow_\cdot(X,B)& \stackrel{\eta^*}{\lraw}& \Flow_\cdot(X_e,B)(-1)&
\lraw& 0\end{array}$$
with exact rows, in which the bottom row is split by $\pi^*$ (and
$\eta^*$ is an abbreviation for $\eta_{X,e}^*$).
The map $\phi''$ is just the restriction of $\phi$, and $\phi'(y_c)
:=\alpha_c$ for all $c\in\EE'$.
The map $\tilde{\iota}$ is the natural inclusion, and to describe
$\tilde{\eta}$ we first define a $B$-algebra isomorphism $\tilde{p}:
B_![\xx]\goesto B_![\yy]$ by 
$$\tilde{p}(x_c^{<j>}):=\left\{\begin{array}{ll}
y_{e'}^{<j>} &\mathrm{if}\ c=e,\\
(y_c+s_cy_{e'})^{<j>} &\mathrm{if}\ c\in\EE\drop\{e\},\end{array}\right.$$
for all $c\in\EE$ and $j\in\NN$;  here, as in Lemma $4.3$(c), $s_c:=\beta_c(X_{e'})$.
Now, with $\tilde{\eta}:=\tilde{p}\circ\partial/\partial x_e$ it is
clear that the top row is exact.  We also define $\tilde{\pi}:B_![\yy](-1)
\goesto B_![\xx]$ by putting $\tilde{\pi}(P(\yy)):=
\int \tilde{p}^{-1}(P(\yy)) dx_e$ for all $P(\yy)\in B_![\yy]$.
From the definitions it is clear that $(\tilde{\pi}\circ\tilde{p})(P(\xx))=
\int P(\xx) dx_e$ for every $P(\xx)\in B_![\xx]$.  It is also clear that
$\tilde{\eta}\circ \tilde{\pi}=1$, and so the top row of
the diagram is also split.  It remains to check commutativity of the diagram,
which we do in the following three claims.

\emph{Claim $1$}:\ \ $\phi'\circ\tilde{p}=p_e\circ\phi$.  Since these maps
are $B$-algebra homomorphisms, the following verification suffices.  For
each $c\in\EE$ and $j\in\NN$, using Lemma $4.3$(c), we have
\showonq
(p_e\circ\phi)(x_c^{<j>})
&=& p_e(\beta_c^{<j>}) = \left\{\begin{array}{ll}
\alpha_{e'}^{<j>} & \mathrm{if}\ c=e',\\
(\alpha_c+s_c\alpha_{e'})^{<j>} & \mathrm{if}\ c\neq e',\\
\end{array}\right.\\
&=& \left\{\begin{array}{ll}
\phi'(y_{e'}^{<j>}) & \mathrm{if}\ c=e',\\
\phi'((y_c+s_cy_{e'})^{<j>}) & \mathrm{if}\ c\neq e',\\
\end{array}\right.\\
&=& (\phi'\circ\tilde{p})(x_c^{<j>}).
\showoffq

\emph{Claim $2$}:\ \ $\phi\circ\tilde{\pi}=\pi^*\circ\phi'$.
Consider any $Q(\yy)\in B_![\yy]$.  Since $\tilde{p}$ is surjective, 
there is an $S(\xx)\in B_![\xx]$ such that $\tilde{p}(S(\xx))=Q(\yy)$.
Since the maps are $B$-linear it suffices to verify that
$(\phi\circ\tilde{\pi}\circ\tilde{p})(\xx^{<\jj>})=
(\pi^*\circ\phi'\circ\tilde{p})(\xx^{<\jj>})$ for all $\jj:\EE\goesto\NN$.
Now
$$(\phi\circ\tilde{\pi}\circ\tilde{p})(\xx^{<\jj>})=
\phi\left(\int\xx^{<\jj>} dx_e\right)=\BETA^{<\jj+\delta_e>}$$
and
\showonq
& & (\pi^*\circ\phi'\circ\tilde{p})(\xx^{<\jj>})=
(\pi^*\circ p_e\circ\phi)(\xx^{<\jj>})\\
&=& (\pi^*\circ p_e)(\BETA^{<\jj>})=\BETA^{<\jj+\delta_e>},
\showoffq 
using Claim $1$ and Lemma $4.3$(d).

\emph{Claim $3$}:\ \ $\phi'\circ\tilde{\eta}=\eta^*\circ\phi$.
Since the maps are $B$-linear it suffices to verify the identity
applied to each $\xx^{<\jj>}$, for $\jj:\EE\goesto\NN$, so fix any
such $\jj$.  First notice that
$$(\tilde{\pi}\circ\tilde{\eta})(\xx^{<\jj>})=\int\left(\frac{\partial
\xx^{<\jj>}}{\partial x_e}\right) dx_e = \tv[j(e)>0]\xx^{<\jj>}.$$
Now, since $\eta^*\circ\pi^*=1$ on $\Flow_\cdot(X_e,B)(-1)$,
and since $\eta^*(\BETA^{<\jj>})=0$ if $j(e)=0$, we calculate (using Claim $2$)
that
\showonq
(\phi'\circ\tilde{\eta})(\xx^{<\jj>})
&=& (\eta^*\circ\pi^*\circ\phi'\circ\tilde{\eta})(\xx^{<\jj>})\\
&=& (\eta^*\circ\phi\circ\tilde{\pi}\circ\tilde{\eta})(\xx^{<\jj>})\\
&=& (\eta^*\circ\phi)(\tv[j(e)>0]\xx^{<\jj>})\\
&=& (\eta^*\circ\phi)(\xx^{<\jj>}),
\showoffq
verifying Claim 3.

Since $\phi\circ\tilde{\iota}=i_*\circ\phi''$ is clear, we see that
the diagram is commutative, and since $\phi''$ is surjective it
follows from the kernel-cokernel exact sequence and the induction
hypothesis that
$$\ker(\phi)=\tilde{\iota}(I(X\drop\{e\},Y,\orn))\oplus\tilde{\pi}
(I(X_e,Y',\orn)(-1)).$$
It is clear from the definition that $I(X,Y,\orn)\subseteq\ker(\phi)$.
For the converse, let $Q(\xx)$ be an arbitrary element of $\ker(\phi)$.
From the above, we may write 
$$Q(\xx)=\tilde{\iota}(Q''(\xx))+\tilde{\pi}(Q'(\yy))$$
with $Q''(\xx)\in I(X\drop\{e\},Y,\orn)$ and
$Q'(\yy)\in I(X_e,Y',\orn)(-1)$.  The first term is easily dealt
with, since $\tilde{\iota}$ includes $I(X\drop\{e\},Y,\orn)$ into
$I(X,Y,\orn)$.  For the second term, $Q'(\yy)$ may be
written as a $B$-linear combination of elements of the form
$\yy^{<\jj>} P_{\theta'}^{<1+r>}(\yy)$ with $\jj:\EE'\goesto\NN$,
$\theta'\in\mathcal{C}(X_e)$, and $r\geq\#\supp(\theta')$, so it
suffices to consider just one term of this form.
Since $p_e$ is surjective, there is some $\theta\in\Flow_1(X,B)$
such that $p_e(\theta)=\theta'$.  In fact, since $\theta(X_c)=
\theta'((X_e)_c)$ for all $c\in E(X)\drop\{e\}$, and then $\theta(X_e)$ is
determined by the condition that $\theta$ is a $B$-flow in $X$, this $\theta$
is unique.  Moreover, it follows that since $\theta'\in\mathcal{C}(X_e)$
we have $\theta\in\mathcal{C}(X)$, that $\#\supp(\theta')\leq\#\supp(\theta)
\leq\#\supp(\theta')+1$, and that $\tilde{p}(P_\theta(\xx))=P_{\theta'}(\yy)$.
Now let $L(\xx):= \tilde{p}^{-1}(\yy^{<\jj>})$.  If $\theta(X_e)=0$ then
$\#\supp(\theta)=\#\supp(\theta')$ and
\showonq
\tilde{\pi}(\yy^{<\jj>} P_{\theta'}^{<1+r>}(\yy))
&=& (\tilde{\pi}\circ\tilde{p})(L(\xx)P_\theta^{<1+r>}(\xx))\\
&=& \int L(\xx) P_\theta^{<1+r>}(\xx) dx_e\\
&=& P_\theta^{<1+r>}(\xx)\int L(\xx) dx_e\in I(X,Y,\orn),
\showoffq
since $B_![\xx]$ is closed under $\int dx_e$.
If $\theta(X_e)=b\in\{-1,1\}$ then $\#\supp(\theta)=\#\supp(\theta')+1$,
and using integration by parts repeatedly we find that
\showonq
\tilde{\pi}(\yy^{<\jj>} P_{\theta'}^{<1+r>}(\yy))
&=& (\tilde{\pi}\circ\tilde{p})(L(\xx)P_\theta^{<1+r>}(\xx))
= \int L(\xx) P_\theta^{<1+r>}(\xx) dx_e\\
&=& L(\xx)P_\theta^{<2+r>}(\xx)b^{-1}-
\int\left(\frac{\partial L(\xx)}{\partial x_e}\right)
P_\theta^{<2+r>}(\xx)b^{-1}dx_e\\
&=& L(\xx)P_\theta^{<2+r>}(\xx)b^{-1}-
\left(\frac{\partial L(\xx)}{\partial x_e}\right)
P_\theta^{<3+r>}(\xx)b^{-2}\\
& & +\int\left(\frac{\partial^2L(\xx)}{\partial^2x_e}\right)
P_\theta^{<3+r>}(\xx)b^{-2}dx_e\\
&=& \ldots\ =\sum_{i=2}^g M_i(\xx)P_\theta^{<i+r>}(\xx),
\showoffq
for finitely many polynomials $M_i(\xx)\in B_![\xx]$.  Since
$P_\theta^{<i+r>}(\xx)$ is in $I(X,Y,\orn)$ for all $i\geq 2$,
this suffices to complete the proof.
\end{proof}
Notice that if $B$ is a $\QQ$-algebra then $B_![\xx]=B[\xx]$, so that
Theorem $4.8$ presents $\Flow_\cdot(X,B)$ as a quotient of a polynomial
algebra in this case.  For example,
$\Flow_\cdot(\mathsf{K}_4,\QQ)$ is isomorphic
to the quotient of $\QQ[x,y,z]$ by the ideal
$$(x^4,y^4,z^4,(x+y)^5,(x+z)^5,(y+z)^5,(x+y+z)^4)$$
(and has Poincar\'e polynomial $D_{\mathsf{K}_4}(t)=
1+3t+6t^2+10t^3+11t^4+6t^5+t^6$).

\begin{CORO} Let $X$ be a graph with $n$ vertices, $m$ edges, $\ell$
cut-edges, $k$ connected components, and with shortest cycle of length
$g$, and let $D_X(t)=d_0+d_1t+\cdots+d_{m-\ell}t^{m-\ell}$.  Then for all
$0\leq j\leq g$,
$$d_j=\binom{m-n+k+j-1}{j}.$$
\end{CORO}
\begin{proof}
Adopt the notation (*) with $B=\QQ$.  The ideal $I(X,Y,\orn)$
of $\QQ[\xx]$ is zero in the graded components of degree $0\leq j\leq g$,
and so by Theorem $4.8$, $d_j$ is the dimension of the $j$-th graded
component of $\QQ[\xx]$ for all $0\leq j\leq g$, proving the result.
\end{proof}

\begin{THM}  Let $X$ be a graph with $m$ edges and $\ell$ cut-edges, and
let $B$ be a commutative ring.  Coordinatize $\Flow_\cdot(X,B)$ by
some $\orn\in\Orn$.  Let $\varphi\in\Flow_1(X,B)$ be such that $\varphi
(X_e)$ is a non-zerodivisor in $B$ for each non-cut-edge $e$ of $X$.
For any integer $0\leq j\leq (m-\ell)/2$, if $\cha(B)$ is zero or coprime to
$\binom{m-\ell-j-i}{j-i}$ for all $0\leq i\leq j$ then the homomorphism
$\Flow_j(X,B)\goesto \Flow_{m-\ell-j}(X,B)$ defined by
$\theta\mapsto\theta\cdot\varphi^{<m-\ell-2j>}$ is injective.
\end{THM}
\begin{proof}  Because of (the dual of) Theorem $1.2$ and the fact that 
$\Flow_\cdot(\cutedge,B)=B$, it suffices to consider the case of a 
graph $X$ with no cut-edges.  Consider an integer $0\leq j\leq m/2$.
Let $M$ be the matrix with rows indexed by
$\Po_{m-j}(X)$ and columns indexed by $\Po_j(X)$, with $M_{\sigma,\rho}:=1$
if $\rho\subseteq\sigma$ and $M_{\sigma,\rho}:=0$ otherwise.
Let $P$ be the matrix which represents the homomorphism $\cdot
\varphi^{<m-2j>}:\Hom(F^j(X),B)\goesto\Hom(F^{m-j}(X),B)$ with
respect to the bases $\{\xi_\rho:\ \rho\in\Po_j(X)\}$ and
$\{\xi_\sigma:\ \sigma\in\Po_{m-j}(X)\}$.
That is, $P_{\sigma,\rho}:=\prod_{e\in\sigma\drop\rho}\varphi(X_e)$ if
$\rho\subseteq\sigma$, and $P_{\sigma,\rho}:=0$ otherwise.
Let $R$ be the diagonal matrix with entries $R_{\rho,\rho}:=\prod_{e\in\rho}
\varphi(X_e)$, and let $S$ be the diagonal matrix with entries
$S_{\sigma,\sigma}:=\prod_{e\in\sigma}\varphi(X_e)$.  One verifies that
the matrix equation $PR=SM$ holds, and hence $\det(P)\det(R)=\det(S)\det(M)$.
By the hypothesis on $\varphi$, both $\det(R)$ and $\det(S)$ are
non-zerodivisors in $B$, and hence $\det(P)$ is a non-zerodivisor if
and only if $\det(M)$ is.   Thus, if $\det(M)$ is coprime to $\cha(B)$
then $\cdot\varphi^{<m-2j>}$ is injective, and so its restriction to
$\Phi_j(X,B)$ is also injective.  From a theorem of Wilson \cite{Wi}
it follows that
$$\det(M)=\prod_{i=0}^j\binom{m-j-i}{j-i}^{\binom mi -\binom m{i-1}}$$
which suffices to complete the proof.  \end{proof}

\begin{CORO}  Let $X$ be a graph with $m$ edges and $\ell$ cut-edges,
and let $D_X(t)= d_0+d_1t+\cdots+d_{m-\ell}t^{m-\ell}$.  Then
$1=d_0\leq d_1\leq\cdots\leq d_{\lfloor (m-\ell)/2\rfloor}$ and
$d_j\leq d_{m-\ell-j}$ for each $0\leq j\leq (m-\ell)/2$. \end{CORO}
\begin{proof}  Apply Theorem $4.10$ with $B=\QQ$.  With the notation (*),
if $\mathcal{E}= \{e_1,...,e_d\}$ then $\varphi:=\sum_{i=1}^d 3^i
\beta_{e_i}$ is a $\ZZ$-flow in $X$ satisfying the hypothesis of Theorem
$4.10$.  The inequalities $d_j\leq d_{m-\ell-j}$ are immediate.
Since $\varphi^{<m-2j>}=\varphi^{m-2j}/(m-2j)!$ for all $0\leq j\leq m/2$,
it follows
that $\cdot\varphi:\Flow_j(X,B)\goesto\Flow_{j+1}(X,B)$ is injective for
all $0\leq j\leq (m-\ell)/2$, implying the remaining inequalities.
\end{proof}

\setcounter{section}{4}

\section{Integer flows and theta functions.}

Finally, we determine the determinant and theta function of the lattice
$\Flow_1(X,\ZZ)$ of integer-valued flows on the graph $X$ in terms of
combinatorial data of $X$.  Other results about this lattice
(including the determinant) have been obtained by Bacher, de la Harpe,
and Nagnibeda \cite{BHN};  in particular, it is interesting to compare
their description of the facets of the Voronoi polytope of
$\Flow_1(X,\ZZ)$ with
the generators of the ideal $I(X,Y,\orn)$ described in Section $4$.

For a graph $X$, define a nondegenerate symmetric positive definite
bilinear form $\la\cdot,\cdot\ra$ on $\Hom(F^\cdot(X),\RR)$ by requiring
that the basis $\{\xi_\sigma:\ \sigma\in\Po(X)\}$ be orthonormal.
That is,
$$\la\varphi,\theta\ra:=\sum_{\sigma\in\Po(X)}\varphi(X_\sigma)
\theta(X_\sigma)$$
for all $\varphi,\theta\in\Hom(F^\cdot(X),\RR)$.
This makes $\Hom(F^\cdot(X),\RR)$ into a Euclidean space.
Since $\la\cdot,\cdot\ra$ is invariant under the automorphisms $\Psi$,
it induces a well-defined Euclidean form on $\Flow_\cdot(X,\RR)$,
also denoted by $\la\cdot,\cdot\ra$.

\begin{LMA} Let $X$ be a graph.\\
\textup{(a)}\ For any nonempty affine subspace $U$ of $\Hom(F^\cdot(X),\RR)$,
there is a unique $\varphi\in U$ with $\la\varphi,\varphi\ra$ minimum.\\
\textup{(b)}\ If $U$ is defined by equations with coefficients in
$\QQ$ (relative to the basis $\{\xi_\sigma:\ \sigma\in\Po(X)\}$) then
the $\varphi$ from part \textup{(a)} is in $\Hom(F^\cdot(X),\QQ)$.\\
\textup{(c)}\ Coordinatize $\Flow_\cdot(X,\RR)$ by some $\orn\in\Orn$.
Let $\sigma\in\Po(X)$ contain no cut-edges of $X$.  Then there is a
unique $\varphi\in \Flow_\cdot(X,\QQ)$ such that $\varphi(X_\sigma)=1$
and $\la\varphi,\varphi\ra$ is minimum subject to this condition.\\
\textup{(d)}\ With notation as in \textup{(c)}, if $\#\sigma=j$ then
$\varphi\in\Flow_j(X,\QQ)$.
\end{LMA}
\begin{proof}
Part (a) is the fact that an affine subspace $U$ of a Euclidean space
contains a unique vector $\varphi$ which is closest to the origin;  indeed,
if $U_0$ is the translate of $U$ passing through $\mathbf{0}$ and $U_0^\bot$
is the orthogonal complement of $U_0$ then $\varphi$ is the unique vector
in $U_0^\bot\cap U$.  If $U$ is defined over $\QQ$ then so is $U_0^\bot$, and
it follows that $\varphi$ has rational coordinates, proving (b).  Part
(c) follows from (a) and (b) since if $\sigma$ contains no cut-edges
of $X$ then $U:=\lambda_\orn(\Flow_\cdot(X,\QQ))\cap\{\varphi\in
\Hom(F^\cdot(X),\QQ):\ \varphi(X_\sigma)=1\}$ is not empty, and part (d)
is evident.
\end{proof}
 
Consider an arc $a$ such that $\und{a}$ is not a cut-edge of $X$, 
coordinatize $\Flow_\cdot(X,\RR)$ by some $\orn\in\Orn$ such that
$s_\orn(a)=1$, and let $\chi^\orn_a$ be the element of $\Flow_1(X,\QQ)$
guaranteed by Lemma 5.1(c,d).  One checks that if $\eps\in\Orn$ is such
that $s_\eps(a)=1$ then $\chi^\eps_a=\chi^\orn_a\circ\psi_\orn^\eps$.
It follows that $\chi_a:=\lambda_\orn^{-1}(\chi^\orn_a)$ does not depend
on the choice of $\orn$;  this is the \emph{characteristic flow of $a$ in
$X$}.  It is easy to see that $\chi_{\rev{a}}=-\chi_a$.

\begin{LMA}  Let $X$ be a graph, let $a\in A$ be such that $e:=\und{a}$
is not a cut-edge, and let $\chi_a$ be the characteristic flow of $a$ in $X$.
Coordinatize $\Flow_\cdot(X,\RR)$ by some $\orn\in\Orn$ such that
$s_\orn(a)=1$.\\
\textup{(a)}\ For any $\varphi\in\Flow_1(X,\RR)$, one has $\la\varphi,\chi_a
\ra=\varphi(X_e)\la\chi_a,\chi_a\ra$.\\
\textup{(b)}\ There is a ``potential function'' $\nu_a:V\goesto\QQ$ such that
$\nu_a(t(a))=0$ and for all $c\in E\drop\{e\}$,
$\chi_a(X_c)=\nu_a(t(\orn(c)))-\nu_a(o(\orn(c)))$.\\
\textup{(c)}\ With $\nu_a$ from part \textup{(b)}, we have
$\la\chi_a,\chi_a\ra=1+\nu_a(o(a))$.
\end{LMA}
\begin{proof}
For part (a), first consider the case that $\varphi(X_e)=0$, and suppose that
$\la\varphi,\chi_a\ra\neq 0$.  Replacing $\varphi$ by $-\varphi$ if necessary,
we may assume that $\la\varphi,\chi_a\ra>0$.  Let $\theta:=\chi_a-\eps\varphi$,
with $\eps>0$ small.  Now
$$\la\theta,\theta\ra=\la\chi_a,\chi_a\ra-2\eps\la\varphi,\chi_a\ra
+\eps^2\la\varphi,\varphi\ra,$$
and as $\eps\goesto 0$ the last term is negligible.  Thus, for some 
$\eps>0$ we have $\la\theta,\theta\ra<\la\chi_a,\chi_a\ra$, but since
$\theta\in\Flow_1(X,\RR)$ and $\theta(X_e)=1$, this contradicts the definition
of $\chi_a$.  Thus, $\la\varphi,\chi_a\ra=0$, as claimed.  For the general
case of part (a), consider any $\varphi\in\Flow_1(X,\RR)$ and let
$\varphi':=\varphi-\varphi(X_e)\chi_a$.  Then $\varphi'(X_e)=0$, so
$\la\varphi',\chi_a\ra=0$, so $\la\varphi,\chi_a\ra=\varphi(X_e)\la\chi_a,\chi_a
\ra$, as claimed.

For part (b), consider any $v\in V$, and let $P:=(a_1,...,a_h)$ be a 
sequence of arcs in $A\drop\{a,\rev{a}\}$ such that $o(a_1)=t(a)$,
$t(a_h)=v$, and $t(a_i)=o(a_{i+1})$ for all $1\leq i\leq h-1$.
Define $\psi_P\in\Hom(F^\cdot(X),\QQ)$ by
$\psi_P:=\sum_{i=1}^h s_\orn(a_i)\xi_{\und{a}_i}$, and let
$\nu_a(v):=\la\psi_P,\chi_a\ra$.  To check that this is well-defined,
let $Q:=(b_1,...,b_g)$ be another walk in $X\drop\{e\}$ from $t(a)$ to $v$,
and define $\psi_Q$ as for $\psi_P$.  Then $\psi_P-\psi_Q$ is in
$\Flow_1(X,\QQ)$ and $(\psi_P-\psi_Q)(X_e)=0$, so $\la\psi_P,\chi_a\ra
=\la\psi_Q,\chi_a\ra$ from part (a), as required.  It is clear that
$\nu_a(t(a))=0$.   Now, for any edge $c\in E\drop\{e\}$,
let $P$ be a walk in $X\drop\{e\}$ from $t(a)$ to $o(\orn(c))$, and let
$Q:=(P,\orn(c))$ be the concatenation of $P$ with $\orn(c)$.  Then
$\chi_a(X_c)=\la(\psi_Q-\psi_P),\chi_a\ra=\nu_a(t(\orn(c)))-\nu_a(o(\orn(c)))$.

For part (c), let $C$ be any cycle in $X$ containing $e$, and let $\theta\in
\Flow_1(C,\QQ)\subseteq\Flow_1(X,\QQ)$ be such that $\theta(X_e)=1$.  Then
$\la\theta,\chi_a\ra=\la\chi_a,\chi_a\ra$, by part (a).  But since the
support of $\theta$ is the cycle $C$ containing $e$, it is easy to see that
$\la\theta,\chi_a\ra=1+\nu_a(o(a))$, from part (b).
\end{proof}
Lemma 5.2 can be interpreted physically by imagining the graph $X$ as an
electrical network in which each edge has unit resistance.  When a unit
current is forced through the edge $e$ in the direction of $a$, the
principle of least action ensures that the current flowing in each edge
of $X$ is given by $\chi_a$.  The electric potential at each vertex is
given by $-\nu_a$.  Part (b) is one formulation of Kirchhoff's Second Law.

For a graph $X$ let $\cplx(X)$ denote the number of maximal forests of $X$,
often called the \emph{complexity} of $X$.
\begin{PROP}  Let $X$ be a graph, let $a\in A$ be such that $e:=\und{a}$
is not a cut-edge, and coordinatize $\Flow_\cdot(X,\RR)$ by some $\orn\in\Orn$
such that $s_\orn(a)=1$.  Then $\la\chi_a,\chi_a\ra=\cplx(X)/
\cplx(X\drop\{e\})$.
\end{PROP}
\begin{proof}
Let $T(X)$ denote the set of maximal forests of $X$, and notice that
$T(X\drop\{e\})\subset T(X)$ since $e$ is not a cut-edge of $X$.  Let
$\QQ T(X)$ denote the $\QQ$-vector space with basis $T(X)$.  Define
$L:\QQ T(X\drop\{e\})\goesto \QQ T(X)$ as follows.  For a given $Y\in
T(X\drop\{e\})$ there is a unique directed path $P=(b_1,...,b_h)$ in $Y$ from
$t(a)$ to $o(a)$.  Define
$$L(Y):= Y+\sum_{i=1}^h s_\orn(b_i)\chi_a(X_{\und{b}_i})(Y\cup\{e\}\drop
\{\und{b}_i\})$$
(and extend the definition $\QQ$-linearly).  Now consider
$$\sum_{Y\in T(X\drop\{e\})} L(Y) =\sum_{Z\in T(X)} c_Z Z,$$
for some coefficients $c_Z\in\QQ$.  If $Z\in T(X\drop\{e\})$ then
clearly $c_Z=1$.  On the other hand, if $e\in E(Z)$ then 
$c_Z=\sum_b s_\orn(b)\chi_a(X_{\und{b}})$, in which the sum is over all arcs
$b\not\in e$ such that $o(b)$ is in the component of $Z\drop\{e\}$ containing
$t(a)$, and $t(b)$ is in the component of $Z\drop\{e\}$ containing $o(a)$.
By Lemma $1.1$, and since $\chi_a$ is a flow, it follows that
$c_Z-\chi_a(X_e)=0$.  Thus we have
$$\sum_{Y\in T(X\drop\{e\})} L(Y) =\sum_{Z\in T(X)}Z.$$
Now apply the linear functional on $\QQ T(X)$ defined by $Z\mapsto 1$ for all
$Z\in T(X)$.  For $Y\in T(X\drop\{e\})$ with path $P=(b_1,...,b_h)$ as above,
Lemma 5.2 implies that
\showonq
L(Y) &\mapsto& 1+\sum_{i=1}^h s_\orn(b_i)\chi_a(X_{\und{b}_i})\\
&=& 1+\sum_{i=1}^h s_\orn(b_i)(\nu_a(t(b_i))-\nu_a(o(b_i)))\\
&=& 1+\nu_a(o(a))=\la\chi_a,\chi_a\ra.
\showoffq
Therefore $\la\chi_a,\chi_a\ra\cplx(X\drop\{e\})=\cplx(X),$
completing the proof.
\end{proof}

An (integral) \emph{lattice} in a Euclidean space $(W,\la\cdot,\cdot\ra)$ is a
subgroup $\Lambda$ of $(W,+)$ such that $\span_\RR\Lambda=W$
and $\la\cdot,\cdot\ra$ is integer-valued when restricted to $\Lambda$;
Conway and Sloane \cite{CS} is the canonical reference.
For an ordered basis $(\lambda_1,...,\lambda_d)$ of $\Lambda$, its
\emph{Gram matrix} is the $d$-by-$d$ matrix $G$ with entries $G_{i,j}
:=\la\lambda_i,\lambda_j\ra$.  The \emph{determinant of $\Lambda$} is
$\det(\Lambda):=\det(G)$, and does not depend on the choice of ordered
basis.  (In fact, $\det(\Lambda)$ is the square of the volume of the
compact group $W/\Lambda$ with respect to the metric induced from
$\la\cdot,\cdot\ra$.) The \emph{theta function} of $\Lambda$ is
$$\thfn(\Lambda|z):= \sum_{\lambda\in\Lambda}
e^{\pi i\la\lambda,\lambda\ra z}
=\sum_{\lambda\in\Lambda}q^{\la\lambda,\lambda\ra},$$
in which $q:=e^{\pi iz}$. 
We also need theta functions of translates of a lattice $\Lambda$ by
a vector $\gamma$ such that $\la\gamma,\lambda\ra\in\QQ$ for all
$\lambda\in\Lambda$, which merely requires the use of rational exponents
of $q$ (with a bounded denominator for any given $\gamma$).  It is
convenient to express the theta function of $\Flow_1(X,\ZZ)$
in terms of the series
$$\thpsi(\alpha|z):=\sum_{n\in\ZZ} q^{(n+\alpha)^2}$$
(which converge absolutely for $\alpha\in\RR$ and $z\in\CC$ with Im$(z)>0$,
uniformly on any compact subset).  In terms of the Jacobi theta functions
$$\vartheta_3(\xi|z):=\sum_{n\in\ZZ}e^{2i\xi n+\pi in^2z}$$
we have $\thpsi(\alpha|z)=e^{\pi i\alpha^2z}\vartheta_3(\pi\alpha z|z).$
The Jacobi triple product formula (see equation (32) in Section $4.4$ of
\cite{CS}, or Theorem $2.8$ of Andrews \cite{An}) yields
$$\thpsi(\alpha|z)=q^{\alpha^2}\prod_{n=1}^\infty
\left(1-q^{2n}\right)\left(1+q^{2n-1+2\alpha}\right)
\left(1+q^{2n-1-2\alpha}\right).$$

For a graph $X$ and edge $e\in E$, we define the \emph{index of $e$ in $X$}
to be the least positive integer $r$ such that $r\chi_a\in\Flow_1(X,\ZZ)$,
in which $\chi_a$ is the characteristic flow of $a\in e$ in $X$ (since
$\chi_{\rev{a}}=-\chi_a$, the choice of $a\in e$ does not matter).

\begin{THM}  Let $X$ be a graph, and adopt the notation \textup{(*)}
with $B=\ZZ$.
Let $\EE=\{e_1,...,e_d\}$, and for $1\leq i\leq d$ let
$X_i:=X\drop\{e_1,...,e_{i-1}\}$ and let $\chi_i$ be the
characteristic flow of $\orn(e_i)$ in $X_i\subseteq X$.\\
\textup{(a)}\ Then $\{\chi_1,...,\chi_d\}$ are pairwise orthogonal in
$\Flow_1(X,\QQ)$.\\
\textup{(b)}\ The determinant of $\Flow_1(X,\ZZ)$ is $\cplx(X)$.\\
\textup{(c)}\ For $1\leq i\leq d$ let $r_i$ be the index of $e_i$
in $X_i$, and let $\varphi_i:=r_i\chi_i$.  Then 
$$\EuScript{S}:=\left\{\sum_{i=1}^d g_i\beta_i:\ 0\leq g_i<r_i,\ g_i\in\ZZ\right\}$$
is a system of representatives for the cosets of the subgroup $\Gamma$ of
$\Flow_1(X,\ZZ)$ generated by $\{\varphi_1,...,\varphi_d\}$.\\
\textup{(d)}\ For $1\leq i\leq d$ let $w_i:=\la\varphi_i,\varphi_i\ra$.
The theta function of $\Phi_1(X,\ZZ)$ is
$$\thfn(\Flow_1(X,\ZZ)|z)=\sum_{\lambda\in\EuScript{S}}\prod_{i=1}^d
\thpsi\left(\left.\frac{\la\lambda,\varphi_i\ra}{w_i}\right| w_iz\right).$$
\end{THM}
\begin{proof}
For $1\leq h<i\leq d$ we have $\chi_i\in\Flow_1(X_h,\ZZ)$ and 
$\chi_i((X_h)_{e_h})=0$, so by Lemma 5.2(a), $\la\chi_i,\chi_h\ra=0$,
proving part (a).

For part (b), let $P$ be the $m$-by-$d$ matrix
in which the $i$-th column is the coordinate vector of $\chi_i$
with respect to the basis $\{\xi_e:\ e\in\Po_1(X)\}$ of $\Hom(F^1(X),
\QQ)$.  Let $Q$ be the $m$-by-$d$ matrix in which the $i$-th column
is the coordinate vector of $\beta_i$ with respect to the same basis.
By Theorem 4.4, $(\beta_1,...,\beta_d)$ is an ordered basis for
$\Flow_1(X,\ZZ)$, with Gram matrix $Q^\top Q$.  For each $1\leq i\leq d$
we have $\chi_i=\sum_{h=1}^d \chi_i(X_{e_h})\beta_h$;  denoting by $M$ the
$d$-by-$d$ matrix such that $M_{h,i}:=\chi_i(X_{e_h})$ we have
$P=QM$.  By construction of $(\chi_1,...,\chi_d)$, $M$ is lower
triangular and $M_{i,i}=1$ for all $1\leq i\leq d$.  Therefore,
$\det(M)=1$ and hence
$$\det(\Flow_1(X,\ZZ)) = \det(Q^\top Q)=\det((M^{-1})^\top P^\top P M^{-1})
= \det(P^\top P).$$
But part (a) shows that $P^\top P$ is diagonal, and by Proposition $5.3$,
$$\det(P^\top P)=\prod_{i=1}^d\frac{\cplx(X_i)}{\cplx(X_{i+1})}=\cplx(X),$$
in which we have defined $X_{d+1}:=Y$, so that $\cplx(X_{d+1})=1$.

For part (c), every $\theta\in\Flow_1(X,\ZZ)$ may be written as
$\theta=\sum_{i=1}^d f_i\beta_i$ for unique integers $f_1,...,f_d$
(in fact, $f_i=\theta(X_{e_i})$).  By the division algorithm we write
$f_1=h_1r_1+g_1$ with $0\leq g_1<r_1$, and consider
$\theta_2:=\theta-h_1\varphi_1$.  Inductively, for $2\leq i\leq d$
we write $f_i=h_ir_i+g_i$ with $0\leq g_i<r_i$, and consider
$\theta_{i+1}:=\theta_i-h_i\varphi_i$.  Since $\varphi_i(X_{e_h})=0$
if $h<i$, it follows that $\theta_{d+1}\in\EuScript{S}$, the candidate set of
representatives.  Since $\theta=h_1\varphi_1+\cdots+h_d\varphi_d+
\theta_{d+1}$, this shows that every $\ZZ$-flow in $X$ is congruent
modulo $\Gamma$ to some flow in $\EuScript{S}$.  Conversely, assume that
$\theta,\theta'\in\EuScript{S}$ are
congruent modulo $\Gamma$, and suppose that $\theta\neq\theta'$.
Let $i$ be the least index such that $\theta(X_{e_i})\neq\theta'
(X_{e_i})$.  From the definition of $\EuScript{S}$ it follows that
$1-r_i\leq (\theta-\theta')(X_{e_i})\leq r_i-1$.  But since
$(\theta-\theta')(X_{e_h})=0$ for all $0\leq h<i$ and $\theta-\theta'$
is assumed to be in $\Gamma$, it follows that $(\theta-\theta')(X_{e_i})$
is a multiple of $r_i$.  Therefore, $\theta(X_{e_i})=\theta'(X_{e_i})$,
a contradiction.  This proves part (c).

For part (d), first notice that for $\alpha\in\QQ$ and $c\in\RR$ such that
$c^2\in\ZZ$, the theta function of $c\alpha+c\ZZ$ is $\thpsi(\alpha|c^2z)$.
Secondly, for translates of lattices
$\lambda_1+\Lambda_1$ and $\lambda_2+\Lambda_2$ in orthogonal
Euclidean spaces, the Pythagorean theorem implies that
$$\thfn((\lambda_1\oplus\lambda_2)+(\Lambda_1\oplus\Lambda_2)|z)
=\thfn(\lambda_1+\Lambda_1|z)\thfn(\lambda_2+\Lambda_2|z).$$
Now, since the cosets of $\Gamma$ partition $\Flow_1(X,\ZZ)$ we have
$$\thfn(\Flow_1(X,\ZZ)|z)=\sum_{\lambda\in\EuScript{S}}\thfn(\lambda+\Gamma|z).$$
Since $\{\varphi_1,...,\varphi_d\}$ are pairwise orthogonal, for each
$\lambda\in S$ we have an orthogonal sum
$$\lambda+\Gamma=\bigoplus_{i=1}^d\left(\frac{\la\lambda,\varphi_i\ra}
{w_i}+\ZZ\right)\varphi_i,$$
since $w_i^{-1}\la\lambda,\varphi_i\ra\varphi_i$ is the projection of
$\lambda$ onto the line spanned by $\varphi_i$.  
From the above remarks it is clear that
$$\thfn(\lambda+\Gamma|z)=\prod_{i=1}^d
\thpsi\left(\left.\frac{\la\lambda,\varphi_i\ra}{w_i}\right| w_i z\right),$$
completing the proof.
\end{proof}
Theorem $5.4$(b) follows from Theorems $2$ and $3$ of Biggs \cite{Bi2},
and also appears as Proposition $1$(iii) of Bacher, de la Harpe, and
Nagnibeda \cite{BHN}.

\setlength{\unitlength}{0.7mm}
\begin{figure}\begin{center}
\begin{picture}(170,90)
\thicklines
\put(20,10){\circle*{3}} \put(110,10){\circle*{3}}
\put(60,10){\circle*{3}} \put(150,10){\circle*{3}}
\put(40,30){\circle*{3}} \put(130,30){\circle*{3}}
\put(20,50){\circle*{3}} \put(110,50){\circle*{3}}
\put(60,50){\circle*{3}} \put(150,50){\circle*{3}}
\put(40,80){\circle*{3}} \put(130,80){\circle*{3}}
\put(20,10){\line(0,1){40}}\put(110,10){\line(0,1){40}}
\put(20,10){\line(1,0){40}}\put(110,10){\line(1,0){40}}
\put(60,50){\line(0,-1){40}}\put(150,50){\line(-1,0){40}}
\put(60,50){\line(-1,0){40}}\put(150,50){\line(0,-1){40}}
\put(20,10){\line(1,1){20}}\put(110,10){\line(1,1){20}}
\put(20,50){\line(1,-1){20}}\put(130,30){\line(1,-1){20}}
\put(40,30){\line(1,1){20}}\put(110,50){\line(2,3){20}}
\put(60,50){\line(-2,3){20}}\put(150,50){\line(-2,3){20}}
\qbezier(20,50)(20,65)(40,80)\qbezier(20,50)(40,65)(40,80)
\qbezier(110,50)(115,35)(130,30)\qbezier(110,50)(125,45)(130,30)
\end{picture}\end{center}
\caption{A pair of nontrivially codichromatic graphs.}
\end{figure}
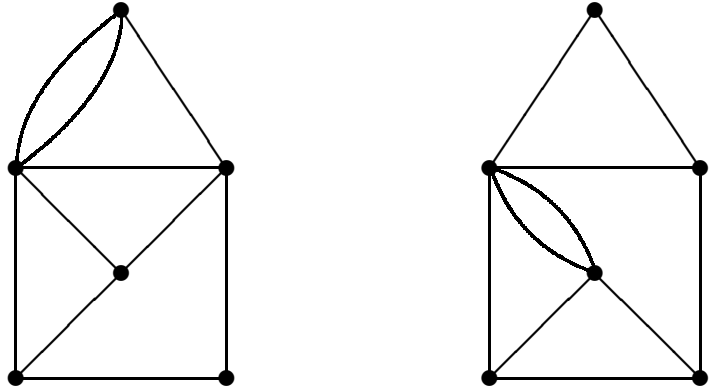

The two graphs depicted in Figure $1$ are \emph{nontrivially
codichromatic}:  they have identical Tutte polynomials but nonisomorphic
graphic matroids.  (This example is due to M.C. Gray; see Tutte \cite{Tu3}.)
Thus, $\Flow_\cdot(X_1,\ZZ)$ and $\Flow_\cdot(X_2,\ZZ)$ are isomorphic
as graded abelian groups.  However, using the Euclidean forms
$\la\cdot,\cdot\ra$ these groups are not isometric.  Indeed, an easy case
analysis shows that $X_1$ (on the left) has $20$ integer-valued flows of
squared-norm $7$ while $X_2$ has $22$ such flows; therefore
$\thfn(\Flow_1(X_1,\ZZ)|z)\neq \thfn(\Flow_1(X_2,\ZZ)|z)$.

\section{Open Problems.}

\begin{PROB}\textup{
Can one identify a basis for $\Flow_j(X,\QQ)$ when $j>1$?  With the notation of
Theorem $4.4$ it would be especially interesting to find a set $\EuScript{J}$
of functions $\jj:\mathcal{E}\goesto\NN$ such that the monomials
$\{\BETA^\jj:\ \jj\in\EuScript{J}\}$ form a basis for $\Flow_\cdot(X,\QQ)$
which is closed under the relation of divisibility.  One hint is
that the dimension of $\Flow_\cdot(X,\QQ)$ as a $\QQ$-vector space is
$D_X(1)=T_X(1,2)$, which is the number of spanning subgraphs of $X$ with the
same number of connected components as $X$ (see, \emph{e.g.} Proposition
$6.2.11$(iii) of Brylawski and Oxley \cite{BrOx}).
If such a basis can be used to obtain results analogous to those in
Section $5$ for $\Flow_j(X,\ZZ)$ with $j>1$ then so much the better.
}\end{PROB}

\begin{PROB}\textup{
For $i,j\geq 0$ with $i+j\leq m-\ell$,
$\Flow_i(X,\ZZ)\cdot\Flow_j(X,\ZZ)$ generates a subgroup of
$\Flow_{i+j}(X,\ZZ)$; denote the quotient by $T_{i,j}(X)$.
Corollary $4.5$ implies that
$T_{i,j}(X)\otimes\QQ=0$, so $T_{i,j}(X)$ is a finite abelian group.
Can one determine the structure of $T_{i,j}(X)$, or even just its order?
(For the $n$-cycle $\Cyc_n$ we have $T_{i,j}(\Cyc_n)=\ZZ/\binom{i+j}{i}\ZZ$
for all $i,j\geq 0$ with $i+j\leq n$.)  Can these groups be used to show
that the ring structure of $\Flow_\cdot(\cdot,\ZZ)$ distinguishes between
pairs of nontrivially codichromatic graphs as in the example at the end
of Section 5?
}\end{PROB}

\begin{PROB}\textup{
Circulation algebras of graphs have some formal similarities with
cohomology rings.  Is there a contravariant functor $X\mapsto\EuScript{X}$
from the category of graphs to the category of topological spaces so that 
$\Flow_\cdot(X,\ZZ)$ is the singular cohomology ring of $\EuScript{X}$?
Alternatively, can $\Flow_\cdot(X,\QQ)$ be interpreted as the rational Chow
ring of some algebraic variety associated to $X$?  (This possibility is
suggested by analogy with the application of Fulton and Sturmfels'
presentation \cite{FS} of the Chow ring of a toric variety to the varieties
associated with partial orders in \cite{W1}.)  It would even be
interesting just to identify an abstract differential graded algebra
with cohomology $\Flow_\cdot(X,\ZZ)$.
}\end{PROB}

\begin{PROB}\textup{
Let $X$ and $Y$ be two graphs with nonempty intersection.  Can one give
an ``excision'' exact sequence (or spectral sequence) for computing
$K^\cdot(X\cup Y)$ in terms of $K^\cdot(X)$, $K^\cdot(Y)$, and
$K^\cdot(X\cap Y)$?  A solution to Problem $6.3$ would provide a big
hint for this.
}\end{PROB}

\begin{PROB}\textup{
Quantum deformation of rational Chow rings of homogeneous spaces is a
hot topic in algebraic geometry these days (see, \emph{e.g.}, Fulton
and Pandharipande \cite{FP}).  Even without a solution to Problem $6.3$,
one could attempt to define a quantum deformation of a circulation
algebra purely combinatorially.  Are there combinatorial
analogues of Gromov-Witten invariants in this context?  If so, then what
do they mean?
}\end{PROB}

\begin{PROB}\textup{
Can one describe the minimal free resolution of $\Flow_\cdot(X,\QQ)$
in terms of the combinatorial structure of $X$?  In particular, is there
a reasonable formula for the dimension of the $j$-th
graded component of the $i$-th resolvent of $\Flow_\cdot(X,\QQ)$?
As a small first step, can one identify a minimal set of generators for
the ideal $I(X,Y,\orn)$ in Theorem $4.8$?
}\end{PROB}

\begin{PROB}\textup{
Can one use circulation algebras to give a proof of the nowhere-zero
$5$-flow conjecture (see Seymour \cite{Se})?  This is rather unlikely,
as a detailed combinatorial argument seems inevitable, but here is a
related question.
Let $\varphi\in\Flow_1(X,\ZZ)$ be such that $\supp(\varphi)=E$ and
$\la\varphi,\varphi\ra$ is minimum with respect to this condition.
Can one give an upper bound on $|\varphi(X_e)|$ independent of $e$ and $X$?
(We have tacitly coordinatized $\Flow_\cdot(X,\ZZ)$ by some orientation.)
In particular, can one prove that $|\varphi(X_e)|<5$?  This would
imply the $5$-flow conjecture.
}\end{PROB}

\begin{PROB}\textup{
How much of the theory of Kirchhoff groups and circulation algebras can
be extended from graphs to more general matroids?  It seems likely that
it will all work for regular (unimodular) matroids, and perhaps even for
all orientable matroids.  For a matroid represented over a field,  
the rowspace of the representing matrix is a natural analogue of $N^1(X)$,
but this depends on the representation.  Can one find for a representable
matroid an analogue of the Kirchhoff group which does not depend on a
choice of representation?  Can anything be said for nonrepresentable matroids?
}\end{PROB}

\begin{PROB}\textup{
The previous problem brings to mind the matroid-theoretic duality between
contraction and deletion.  Contraction is ubiquitous in the above theory,
but deletion is used much less.  Can one make use of deletion to enrich
this theory?   One natural approach is to associate to a graph $X$ not just
$\Flow_\cdot(X,B)$, but the whole family $\{\Flow_\cdot(Y,B):\ Y\subseteq X\}$
as $Y$ ranges over the spanning subgraphs of $X$.  Because of the
algebra homomorphisms induced by graph morphisms, this is a directed set
of algebras and behaves functorially.  Does it have any interesting
or useful structure?
}\end{PROB}

\begin{PROB}\textup{
Many specializations of the Tutte polynomial of a matroid are conjectured
to have a logarithmically concave sequence of coefficients.  Computations
of hundreds of examples suggest that the Poincar\'e polynomials of graphs
also share this property.  That is, for a graph $X$ with $m$ edges and
$\ell$ cut-edges, can one prove that $d_j(X)^2\geq d_{j-1}(X)d_{j+1}(X)$
for all $1\leq j\leq m-\ell-1$?
}\end{PROB}


\begin{thebibliography}{ABC}

\bibitem{An} G.E. Andrews, \emph{The Theory of Partitions}, Encyclopedia 
of Math. and Appl. \textbf{2}, Addison-Wesley, London, Reading MA, 1976.

\bibitem{BHN} R. Bacher, P. de la Harpe, and T. Nagnibeda, \emph{
The lattice of integral flows and the lattice of integral cuts on a 
finite graph}, preprint.

\bibitem{Bi} N. Biggs, \emph{Algebraic Graph Theory}, Cambridge Univ. Press,
Cambridge, 1974.

\bibitem{Bi2} N. Biggs, \emph{Homological coverings of graphs},
J. London Math. Soc. Ser. 2 \textbf{30} (1984), 1-14.

\bibitem{Bj} A. Bj\"orner, \emph{Topological methods}, in: \emph{Handbook
of Combinatorics, vol. II} (R.L. Graham, M. Gr\"otschel, L. Lov\'asz, eds.),
Elsevier Science B.V. and the MIT Press, Amsterdam and Cambridge, 1995.

\bibitem{BCN} A.E. Brouwer, A.M. Cohen, and A. Neumaier, \emph{
Distance Regular Graphs}, Springer-Verlag, 1989.

\bibitem{BrOx} T.H. Brylawski and J.G. Oxley, \emph{The Tutte polynomial
and its applications}, in: \emph{Matroid Applications}, (N. White, ed.),
Cambridge Univ. Press, Cambridge, 1992.

\bibitem{CL} G. Clements and B. Lindstr\"om, \emph{A generalization of a
combinatorial theorem of Macaulay}, J. Combin. Theory \textbf{7} (1969),
230-238.

\bibitem{CS} J.H. Conway and N.J.A. Sloane, \emph{Sphere Packings, Lattices,
and Groups}, Grundlehren der mathematischen Wissenschraften \textbf{290},
Springer-Verlag, Berlin, New York, 1988.

\bibitem{CDS} D.M. Cvetkovi\'c, M. Doob, and H. Sachs, \emph{Spectra of
Graphs}, Academic Press, 1980.

\bibitem{FP} W. Fulton and R. Pandharipande, \emph{Notes on stable maps
and quantum cohomology}, preprint.

\bibitem{FS} W. Fulton and B. Sturmfels, \emph{Intersection theory on
toric varieties}, Topology \textbf{36} (1997), 335-353.

\bibitem{Go} C.D. Godsil, \emph{Algebraic Combinatorics}, Chapman and Hall,
New York and London, 1993.

\bibitem{Lo} L. Lov\'asz, \emph{Kneser's conjecture, chromatic number and
homotopy}, J. Combin. Theory Ser. A \textbf{25} (1978), 319-324.

\bibitem{Mac} S. Mac Lane, \emph{Homology}, Grundlehren der 
mathematischen Wissenschaften \textbf{114}, Springer-Verlag,
Berlin and New York, 1975.

\bibitem{Ox} J.G. Oxley, \emph{Graphs and series-parallel networks},
in: \emph{Theory of Matroids}, (N. White, ed.),  Cambridge Univ.
Press, Cambridge, 1986.

\bibitem{Se} P.D. Seymour, \emph{Nowhere-zero flows}, in: \emph{Handbook
of Combinatorics, vol. I} (R.L. Graham, M. Gr\"otschel, L. Lov\'asz, eds.),
Elsevier Science B.V. and the MIT Press, Amsterdam and Cambridge, 1995.

\bibitem{St} R.P. Stanley, \emph{Combinatorics and Commutative Algebra,
Second Edition}, Birkh\"auser, Berlin, Boston, 1996.

\bibitem{Tu} W.T. Tutte, \emph{A ring in graph theory}, Proc. Camb. Phil.
Soc. \textbf{43} (1947), 26-40.

\bibitem{Tu2} W.T. Tutte, \emph{A contribution to the theory of chromatic
polynomials}, Canad. J. Math. \textbf{6} (1954), 80-91.

\bibitem{Tu3} W.T. Tutte, \emph{Codichromatic graphs}, J. Combin. Theory
Ser. B \textbf{16} (1974), 168-174.

\bibitem{W1} D.G. Wagner, \emph{Singularities of toric varieties
associated with finite distributive lattices}, J. Algebraic Combin.
\textbf{5} (1996), 149-165.

\bibitem{W2} D.G. Wagner, \emph{The Tutte dichromate and Whitney homology
of matroids}, preprint available at \texttt{http://math.uwaterloo.ca/~dgwagner}.

\bibitem{Wa} J.W. Walker, \emph{From graphs to ortholattices and equivariant
maps}, J. Combin. Theory Ser. B \textbf{35} (1983), 171-192.

\bibitem{Wi} R.M. Wilson, \emph{A diagonal form for the incidence matrices
of $t$-subsets vs. $k$-subsets}, Europ. J. Combin. \textbf{11} (1990), 609-615.

\bibitem{Wh} H. Whitney, \emph{$2$-isomorphic graphs}, Amer. J. Math.
\textbf{55} (1933), 245-254.

\end{thebibliography}
\end{document}